\documentclass[12pt]{iopart}

\usepackage{iopams}  

\usepackage{graphicx}

\graphicspath{{figures/}}
\usepackage{subfig}
\usepackage[export]{adjustbox}
\usepackage{url}
\usepackage[hidelinks]{hyperref}
\usepackage{indentfirst}
\usepackage{placeins}
\usepackage[dvipsnames]{xcolor}

\expandafter\let\csname equation*\endcsname\relax

\expandafter\let\csname endequation*\endcsname\relax

\usepackage{amsmath}

\newcommand{\argmin}{\mathop{\mathrm{arg\,min}}}

\begin{document}

\title[On Krylov Methods for Large Scale CBCT Reconstruction]{On Krylov Methods for Large Scale CBCT Reconstruction}

\author{Malena Sabaté Landman\textsuperscript{1}\textsuperscript{*}, Ander Biguri\textsuperscript{1}\textsuperscript{*}, Sepideh Hatamikia\textsuperscript{2,3}, Richard Boardman\textsuperscript{4}, John Aston\textsuperscript{5}, Carola-Bibiane Schönlieb\textsuperscript{1}}

\address{\textsuperscript{1}Department of Applied Mathematics and Theoretical Physics (DAMTP), University of Cambridge, Cambridge, UK}

\address{\textsuperscript{2}Research center for Medical Image Analysis and Artificial Intelligence (MIAAI), Department of Medicine, Danube Private University, Krems, Austria}
\address{\textsuperscript{3}Austrian Center for Medical Innovation and Technology (ACMIT), Wiener Neustadt, Austria} 
\address{\textsuperscript{4}\textmu-Vis X-ray Imaging Laboratory, University of Southampton, Southampton, UK}
\address{\textsuperscript{5}Department of Pure Mathematics and Mathematical Statistics (DPMMS), University of Cambridge, Cambridge, UK}

\address{\textsuperscript{*}These authors contributed equally to this work}
\ead{m.sabate.landman@gmail.com and ander.biguri@gmail.com}

% for internal contacting purposes
% Sepideh email: sepidehhatamikia@yahoo.com 
% Richard email: rpb@soton.ac.uk 
\vspace{10pt}
\begin{indented}
\item[]December 2022
\end{indented}

\begin{abstract}
Krylov subspace methods are a powerful family of iterative solvers for linear systems of equations, which are commonly used for inverse problems due to their intrinsic regularization properties. Moreover, these methods are naturally suited to solve large-scale problems, as they only require matrix-vector products with the system matrix (and its adjoint) to compute approximate solutions, and they display a very fast convergence. 
Even if this class of methods has been widely researched and studied in the numerical linear algebra community, its use in applied medical physics and applied engineering is still very limited. e.g. in realistic large-scale Computed Tomography (CT) problems, and more specifically in Cone Beam CT (CBCT). This work attempts to breach this gap by providing a general framework for the most relevant Krylov subspace methods applied to 3D CT problems, including the most well-known Krylov solvers for non-square systems (CGLS, LSQR, LSMR), possibly in combination with Tikhonov regularization, and methods that incorporate total variation (TV) regularization. This is provided within an open source framework: the Tomographic Iterative GPU-based Reconstruction (TIGRE) toolbox, with the idea of promoting accessibility and reproducibility of the results for the algorithms presented. Finally, numerical results in synthetic and real-world 3D CT applications (medical CBCT and \textmu-CT datasets) are provided to showcase and compare the different Krylov subspace methods presented in the paper, as well as their suitability for different kinds of problems.

\end{abstract}

%\maketitle

\section{Introduction}

Computed Tomography (CT) is a very popular imaging technique widely used in medical and scientific applications. In particular, Cone Beam CT (CBCT) has gained significant attention in the last decade, both for medicine, when low dose image guidance is required (e.g. dental imaging, image guided radiation therapy, image guided surgery), but also in scientific applications involving \textmu-CT,  where higher doses are tolerated in favour of a better image reconstruction quality.  {Moreover, since many clinical applications require producing reliable images in real or near real time \cite{leeser2014fast}\cite{HANSEN}, there is a true need for faster available reconstruction methods. This is crucial in CBCT imaging during surgical procedures, where the long time required by most standard algorithms makes their use unfeasible in a standard clinical workflow.}  {%The computational time need for fast CBCT imaging is becoming increasingly challenging particularly in applications that require imaging immediately before or during surgical procedures. 
For example, this is the case in needle-based procedures, where fast CBCT imaging has the potential to accurately image intraoperative anatomy in close proximity to the needle \cite{diagnostics10121068}\cite{Kickuth1} allowing for immediate adjustment in case of misplacement. Other examples can be found in image guided radiotherapy and online radiotherapy, and particularly in particle radiotherapy, where a CBCT image is taken on-site mere seconds before the radiation dose is delivered  \cite{MACKAY2018293}, with a very limited window for both reconstruction and radiation dose planning. %, making the CBCT acquisition and reconstruction times extremely time-constrained.
Finally, optimization of source-detector CBCT trajectories has recently shown great promise in interventional radiology, offering a variety of benefits, including image quality improvement, FOV expansion, radiation dose reduction, metal artifact reduction, and 3D imaging under kinematic constraints. This optimization process is highly dependent on the image reconstruction speed, so the clinical implementation of such methods can only be realized with the use of fast CBCT reconstruction techniques \cite{10.1371/journal.pone.0245508}\cite{Thies1}.} %\ab{Another typical clinical scenario where computational time is crucial is in image guided radiotherapy or in online radiotherapy. In these procedures, particularly for the more precise particle radiotherapy, a CBCT image is taken on site mere seconds before the radiation dose is delivered \cite{MACKAY2018293}, with a very limited window for both reconstruction and radiation dose planning, making the CBCT acquisition and reconstruction times extremely time-constrained.}

 {In order to perform the CT reconstructions of an image from its measured x-ray projections, one needs to study and understand the properties its underlying numerical model.} Mathematically, this can be formulated as finding a solution of a large-scale linear system of the form
\begin{equation} \label{eq:original_pbm}
    A x + e = b, 
\end{equation}
where $A \in \mathbb{R}^{N \times M}$ is the system matrix describing the measurement process, $b \in \mathbb{R}^{N}$ is the vector of measurements and $e \in \mathbb{R}^{N}$ is the modelled additive noise. Note that $N$ is the number of detector pixels multiplied by the number of projection angles, while $M$ is the number of voxels in the image $x$. For more information see, e.g., \cite{ kak2001principles}\cite{doi:10.1137/1.9780898718836} and references therein. There are two main factors that make problem (\ref{eq:original_pbm}) very challenging to solve in practice. First, the problem is ill-posed, i.e. the matrix $A$ has singular values that decay and cluster at zero, without an evident gap between two consecutive values. This means that the recovered solution is very sensitive with respect to small perturbations (e.g. noise) in the measurements, and therefore some regularization (replacing the original problem by a related more stable problem), needs to be applied to obtain a meaningful reconstruction. The methods described in this paper provide different forms of regularization that will be explained and compared in the following sections. Second, in real-world CT applications, equation (\ref{eq:original_pbm}) can be very large-scale, so it is unfeasible to work directly with the matrix $A$ or, in most cases, even construct it and store it. \\ 

In practical CT applications, an approximation of the solution of (\ref{eq:original_pbm})  {is frequently computed using a direct method commonly known as Filtered Backprojection (FBP),} or as FDK in CBCT problems, and named after Feldkamp, Davis and Kress \cite{feldkamp1984practical}.  {However, these algorithms tend to produce image artifacts due to noise amplification related to the ill-posedness of the problem and mismatches between the idealistic models for the x-ray behaviour and the real measurement sampling process}. An alternative to solve problem (\ref{eq:original_pbm}) is to use iterative methods that rely only on matrix-vector products with $A$ and $A^T$ to handle the large-scale nature of these matrices; hence these are also known as matrix-free methods. These techniques have shown to produce reconstructions of better quality  \cite{kataria2018assessment}\cite{mao2019evaluation}\cite{desai2012impact}, particularly in the cases where there are less measurements, or they are noisier. This is especially relevant in medical applications, where reduced measurements lower the amount of damaging x-ray radiation that is given to the patient. \\

This work focuses on Krylov subspace methods, a family of matrix-free algorithms that are very well-known and studied in the numerical linear algebra community but that have found limited use on real-world CT applications so far. This class of methods was first introduced in the 1950s \cite{Hestenes1952MethodsOC}, but it is recently getting very popular for solving inverse inverse problems \cite{https://doi.org/10.48550/arxiv.2105.07221}\cite{10.1002/gamm.202000017}. Conjugate Gradient Least Squares (CGLS) is the most commonly used Krylov method in applied x-ray CBCT, see, for example  \cite{dabravolski2014dynamic}\cite{pengpen2015motion}\cite{lohvithee2021ant}. Moreover, it is sometimes also found in combination with Tikhonov regularization, or within more complex minimization schemes tackling different variational regularizers \cite{kazantsev2015employing}. Mathematically equivalent to CGLS, the more stable algorithm LSQR, has also been used for CBCT \cite{chillaron2020evaluation}\cite{flores2016application}, but is by far less popular than CGLS. Other minimal residual Krylov solvers, such as the generalized minimal residual algorithms (GMRES), have been seldom used in CT \cite{coban2014regularised}. In particular, we want to include in our comparisons recent developments building up from this algorithm, namely ABBA-GMRES \cite{HANSEN2022114352}\cite{10.1117/12.2646511}, which support unmatched backprojectors: a very common problem in large-scale CT. In this same work LSMR is also used as a comparison. Recent developments in hybrid Krylov methods incorporating Tikhonov regularization and total variation regularization have not been used in real-world CT applications to the best of our knowledge. \\

In the following sections the most relevant Krylov subspace methods for 3D CT problems are presented including the most well-known Krylov solvers for non-square systems (CGLS, LSQR, LSMR), possibly in combination with Tikhonov regularization, and recently developed methods that incorporate Total Variation (TV) regularization. Some of these state-of-the art Krylov methods are newly used in real CT applications, and newly used in a 3D setting. The relevant codes are also provided within an open source framework: the Tomographic Iterative GPU-based Reconstruction (TIGRE) toolbox \cite{tigre}, along with reproducible numerical experiments in synthetic and real-world 3D CT applications (for medical CBCT and \textmu-CT datasets) that showcase and compare the different methods presented in the paper. {Note that Section 2 (Methods) recalls the mathematical framework for the methods described in the paper. Reading this section is not necessary to use the algorithms as provided in the toolbox, so the authors suggest to anyone interested only in the direct applications of such methods to skip this section. Examples of results under different CT acquisition modes and samples are given in Section 3 (Results), some general guidance on the use of the different algorithms is provided in Section 4 (Discussion) and a list of the algorithms in the toolbox is provided in Section 5 (Conclusions).}
\section{Methods}
Krylov methods are projection methods, i.e. iterative methods that, at each iteration $k$, are defined to find the best solution $x_k$ belonging to a Krylov subspace of increasing dimension according to different optimality criteria that define each particular Krylov solver. In particular, Krylov subspaces are generated by the linear combination of the first $k$-$1$ powers of a matrix acting on a vector. In this paper, we mainly focus on subspaces where $A^T A$ acts on $A^T b$, denoted as
\begin{equation}
    \mathcal{K}_k (A^T A,A^T  b) = \text{span} \{ A^T b, (A^T A) A^T b, ..., (A^T A)^{k-1} A^T  b\},
\end{equation}
or in variations thereof. Unless explicitly stated otherwise, the computational cost of the presented Krylov methods is dominated by a matrix vector product by $A$ and $A^T$ per iteration. 
\subsection{Least squares problems}\label{sec_lsqr} Note first that (\ref{eq:original_pbm}) might not be consistent,  {i.e. there might not exist a solution $x^*$ such that $Ax^*=b$}, mainly due to the presence of the noise $e$, but also due to small differences between the discretized model and the true underlying physical model governing the measurement process. Therefore, we consider instead the following least squares problem
\begin{equation} \label{eq:lest_squares}
    \hat{x}= \argmin_x \|A x - b\|_2. 
\end{equation}
Note that solving (\ref{eq:lest_squares}) is equivalent to finding the best linear unbiased estimator for the solution assuming $e$ is white Gaussian noise, see e.g. \cite{vogel}. Generally, the noise in CT problems is related to errors of photon counts and modelled with a Poisson distribution: in the large photon counts regime corresponding to most CT problems, it is reasonable to approximate it by a Gaussian. However, problem (\ref{eq:lest_squares}) can still be very sensitive to small perturbation in the measurements, so the solution of (\ref{eq:lest_squares}) might still be a bad reconstruction of the original image. Krylov methods have inherent regularization properties when combined with early stopping, displaying a phenomena called semiconvergence, i.e. the relative error norm of the solution decreases on the first iterations but starts increasing again after the optimal stopping iteration, see e.g. \cite[Section 6.3]{doi:10.1137/1.9780898718836}. In the following subsections the most applicable Krylov methods to solve problem (\ref{eq:lest_squares}) in the context of CT reconstruction are described. 

\subsubsection{CGLS} Conjugate gradient least squares (CGLS), is the most used Krylov method in CT, and it dates back to \cite{Hestenes1952MethodsOC}. It consists on applying the conjugate gradient method to the normal equations associated to (\ref{eq:lest_squares}): $A^T A x = A^T b$. At each iteration $k$, the solution of
\begin{equation} \label{eq:lest_squares_projected}
    x_k= \argmin_{x \in  \mathcal{K}_k (A^T A,A^T  b) } \|A x - b\|_2
\end{equation}
is computed, such that the residual norm $\|r_k\|$, where $r_k = b - A x_k$, decreases monotonically.
\subsubsection{LSQR} The LSQR method is based on the construction of a Krylov subspace using the Golub–Kahan (GK) bidiagonalization process \cite{lsqr}. This process results on a partial decomposition of $A$ of the form $A V_k = U_{k+1} H_k$, where $H_{k} \in \mathbb{R}^{k+1 \times k}$ is bidiagonal, and such that the orthonormal columns of $V_k$ span the Krylov subspace $\mathcal{K}_k (A^T A,A^T  b)$, $U_{k+1}$ has orthogonal columns and $U_{k+1}e_1 = \|b\| e_1$. Then, problem (\ref{eq:lest_squares_projected}) can be reformulated as
\begin{equation}\label{eq:proj_LSQR}
    x_k = V_k y_k \quad \text{where} \quad y_k = \argmin_{y \in \mathbb{R}^k} \|b - A V_k y\|_2 =  \argmin_{y \in \mathbb{R}^k} \|\|b\| e_1 - H_{k} y\|_2,
\end{equation}
where $e_1$ is the canonical vector of appropriate dimension. Even if this method has been less used in applied CT papers compared to CGLS, these two methods are mathematically equivalent, and LSQR was originally designed to provide a more stable algorithm. A detailed implementation of this method based on short recursions can be found in the original paper \cite{lsqr}. 

\subsubsection{LSMR} Similarly to LSQR, LSMR is also based on the construction of a Krylov subspace using the Golub–Kahan (GK) bidiagonalization process \cite{lsmr}. However, at each iteration, LSMR seeks a solution $ x_k \in \mathcal{K}_k (A^T A,A^T  b)$ such that $\|A^T r_k\|$ is minimized, i.e.
\begin{equation} \label{eq:proj_LSMR}
    x_k = V_k y_k \quad \text{where } y_k = \argmin_{y \in \mathbb{R}^k} \|A^T r_k\|_2 =  \argmin_{y \in \mathbb{R}^k} \|\|A^T b\| e_1 - \bar{H}_k^T H_k y\|_2,
\end{equation}
where $\bar{H}_k \in \mathbb{R}^{k \times k}$ corresponds to the first $k$ rows of the matrix $H_k$. Although both LSQR and LSMR converge in exact arithmetic to the same solution, see e.g. \cite{lsmr}, they produce slightly different solutions at each iteration. Moreover, LSMR is mathematically equivalent to GMRES \cite{doi:10.1137/0907058} applied to the normal equations $A^T A x = A^T b$, and since the system matrix for the normal equations is symmetric, this is also equivalent to using MINRES \cite{doi:10.1137/0712047}.

\subsubsection{AB-GMRES and BA-GMRES} Due to how efficient implementations of the CBCT problems are coded for GPUs, the matrix $B$ that represents the backprojection operator, i.e. the adjoint of $A$, is usually just an approximation of $A^Tb$ \cite{BIGURI202052}. This mismatch can cause the standard Simultaneous Iterative Reconstruction Technique (SIRT) family of iterative solvers to diverge, unless specific perturbations are added to stabilize the convergence, see \cite{HANSEN2022114352}. Alternatively, the approximated transpose matrix $B$ can be used as a right (resp. left) preconditioner for GMRES when solving problem (\ref{eq:lest_squares}), giving rise to AB-GMRES (resp. BA-GMRES) \cite{HANSEN2022114352}. 

\subsection{Tikhonov regularization}
Another form of regularization is Tikhonov regularization, and it is perhaps the simplest and most well-known variational regularization method. It consists on computing the solution
\begin{equation}\label{eq:Tikhonov}
    \hat{x}= \argmin_x  \left\{ \|A x - b\|_2^2+\lambda^2 \|x\|^2_2 \right\},
\end{equation}
where $\lambda$ balances the effect of the fit-to-data term $\|A x - b\|_2^2$ (promoting consistency of the solution with the measurements) and the regularization term $ \|x\|^2_2$ (promoting regularity of the solution). If $\lambda$ is chosen adequately, the semiconvergence behaviour can most times be alleviated, and the algorithms are less sensitive to early stopping; moreover, this allows for a bigger Krylov space to be built, sometimes leading to solutions of improved quality with respect to their non-Tikhonov-regularized counterparts. When $\lambda$ is known, or fixed ahead of the iterations, one can apply any iterative solver (e.g. CGLS, LSQR or LSMR) to the augmented system:
\begin{equation}\label{eq:Tikhonov_augmented}
    \hat{x} = \argmin_x \left\| \left[\begin{array}{c}
A \\
\lambda I
\end{array}
\right]x -
\left[\begin{array}{c}
b \\
0
\end{array}
\right] \right\|.
\end{equation}
For example, an LSMR implementation for fixed $\lambda$ is given in the original paper \cite{lsmr} and compared in this study. An alternative approach is to use hybrid methods: which consist on adding Tikhonov regularization to the projected problem (\ref{eq:proj_LSQR}) or (\ref{eq:proj_LSMR}). In the case of LSQR, this is mathematically equivalent to projecting the regularized problem (\ref{eq:Tikhonov_augmented}) \cite[Chapter 6]{doi:10.1137/1.9780898718836}. However, this is not the case for LSMR \cite{hlsmr}. The big advantage of hybrid methods is that they provide a framework to estimate $\lambda$ on-the-fly when it is not known a-priori. Even if they display a very fast convergence, the drawback of these methods is that they come with the additional cost of having to store $k$ additional (basis) vectors for computing the solution at iteration $k$. This makes them suited for small to medium problems, e.g. $x \in \mathbb{R}^{512 \times 512 \times 512}$, $b \in \mathbb{R}^{512 \times 512 \times 360}$. In some cases, this could be alleviated by storing the coefficients and recomputing all the basis vectors at the end of the iterations requiring twice as many matrix-vector products with $A$ and $A^T$ than their non-hybrid counterparts. We provide a version of hybrid LSQR to show the performance of these methods. For more information, a great review on hybrid methods can be found in \cite{https://doi.org/10.48550/arxiv.2105.07221}.

\subsubsection{hybrid LSQR} Using the same Krylov subspace described for LSQR and adding regularization to the projected problem (\ref{eq:proj_LSQR}), leads to solving, at each iteration $k$:
\begin{eqnarray}\label{eq:hybrid_LSQR}
    x_k = V_k y_k \quad \text{where } y_k = \argmin_y \left\{\left\| \|b\| e_1 - H_k y \right\|_2 + \lambda_k^2 \|y\|_2^2 \right\}.
\end{eqnarray} 
As already mentioned, and thanks to the shift invariance property of Krylov subspaces, for fixed $\lambda_k$ problem (\ref{eq:hybrid_LSQR}) is equivalent to projecting problem (\ref{eq:Tikhonov_augmented}) onto the Krylov subspace $\mathcal{K}_k (A^T A, A^T b)$, see, originally \cite{10.1145/355993.356000}, or \cite{doi:10.1137/1.9780898718836} for a more detailed explanation. An interesting feature of formulation (\ref{eq:hybrid_LSQR}) is that $\lambda_k$ can be computed on-the-fly at each iteration $k$ according to a parameter choice criterion; examples of which can be found in the following section. 
\subsubsection{Parameter choice criteria} A good choice of the regularization parameters is crucial to obtain meaningful reconstructions when dealing with ill-posed problems. In this section we focus on choices for $\lambda_k$ (but note that the total amount of iterations $k$ can also be considered a regularization parameter for regularization by early-stopping). In the following, we provide the description of two of the most simple and commonly used regularization parameter choice criteria. This is by no means an exhaustive list of the available options and we point the interested reader to the reviews in e.g. \cite{https://doi.org/10.48550/arxiv.2105.07221}\cite{10.1002/gamm.202000017}. 

If a good estimate of the norm of the error $\|e\|$ is available, a very popular and reliable parameter choice criterion is the  Discrepancy Principle (DP) \cite{Mor}. This method is based on the idea that
\begin{equation}
    \| A x^{exact} - b \|_2^2 = \| A x^{exact} - b^{exact} - e \|_2^2 = \|e\|_2^2 = \frac{\|e\|_2^2}{\|b\|_2^2} \|b\|_2^2 = \text{nl}^2 \, \|b\|_2^2
\end{equation}
so at each iteration, $\lambda_k$ is chosen so that
\begin{equation}
   \lambda_k = \argmin_{\lambda} \{ \| A x_{\lambda} - b \|_2^2 - \text{nl}^2 \, \|b\|_2^2 \}.
\end{equation}

Alternatively, one can use parameter choice rules that do not use any information about the noise $e$, also known as “heuristic methods”. In particular, we provide an implementation of the Generalized Cross Validation (GCV) parameter choice criterion, which relies on cross validation: a well known statistical tool used to predict possible missing data values. In this case, each of the components of the vector of measurements $b$ is estimated using the rest of components, and the regularization parameter $\lambda_k$ associated with the best predicted values is taken at each iteration. In practice, for hybrid LSQR, using GCV involves solving the following minimization:
\begin{equation}
   \lambda_k = \argmin_{\lambda} \frac{\left\|(I - H_k H^{\dagger}_{k,\lambda}) \, \|b\| e_1 \right\|^2}{\text{tr}((I - H_k H^{\dagger}_{k,\lambda}))^2} \quad \text{where } H^{\dagger}_{k,\lambda} =(H_{k}^{T} H_k + \lambda I)^{-1} H_{k}^{T}.
\end{equation}
Note that this can be generalized to other Krylov methods (e.g. LSMR) by replacing the projected matrix and right hand side by corresponding ones (see, e.g. \cite{hlsmr}). 
\subsection{Total Variation (TV) regularization}
Total variation is a very common variational regularization scheme that promotes piecewise-constant reconstructions by favouring solutions with a sparse gradient. This is very popular in imaging problems as it contributes to preserve edges in the reconstructed image. In this paper the discrete isotropic total variation in 3D is considered: 
\begin{equation} \label{eq:TV}
TV(x) = \sum_{i} \sqrt{[D_l x]_i^2+[D_j x]_i^2+[D_k x]_i^2} =  \left\| \sqrt{[D_l x]_i^2+[D_j x]_i^2+[D_k x]_i^2} \right\|_1,
\end{equation}
where $D_l$, $D_j$, $D_k$ refer to the finite difference approximations of the three directional derivatives for the 3D image $x$. A popular approach to solve the TV problem using Krylov methods is to re-write the TV regularization term using a weighted 2-norm:
\begin{equation} \label{eq:TV_reweighted}
    \hat{x} = \argmin_x \left\{ \|A x - b\|_2^2+\lambda^2 TV(x) \right\} =  \argmin_x \left\{ \|A x - b\|_2^2+\lambda^2 \|W(Dx) Dx\|_2\right\},
\end{equation}
where $D$ is the 3D discrete derivative operator and $W(Dx)$ is a (diagonal) weighting matrix that depends on $Dx$. Then, the functional in (\ref{eq:TV_reweighted}) can be approximated locally by a sequence of quadratic functionals, giving rise to a sequence of problems of the form:
\begin{equation} \label{eq:TV_approx}
    x^{(k)} = \argmin_x \left\{ \|A x - b\|_2^2+\lambda^2 \|L^{(k)} Dx\|_2\right\},
\end{equation}
where $L^{(k)}$ are approximations of $W(Dx)$ of improving quality. This scheme is called iteratively reweighted norm (IRN) and was first used in combination with TV in \cite{IRN} for 2D imaging problems. In the following, two algorithms that (partially) solve the problems in (\ref{eq:TV_approx}) to approximate TV regularization are described. 

\subsubsection{CGSL-TV} The sequence of problems (\ref{eq:TV_approx}) can be solved in an inner-outer scheme fashion where, at each outer iteration, the computed solution $x^{(k)}$ is used to update the weights $L^{(k+1)} = W(D x^{(k)})$. Following \cite{IRN}, an adaptation of this method for 3D using CGLS in the inner iterations, is provided in this paper. This scheme has provable convergence guarantees, but requires $\lambda$ to be known a-priori and can be very computationally expensive due to its inner-outer scheme nature. Other variations of this method have been implemented using other Krylov methods for the inner iterations, e.g., in combination with LSQR \cite{IRN_lsqr}. 

\subsubsection{hybrid fLSQR} An equivalent formulation to (\ref{eq:TV_approx}), dropping the $(k)$ upper-script to ease the notation so that $L = L^{(k)}$, is to solve
\begin{equation}\label{eq:TV_standard_form}
\hat{x} = L^{\dagger}_{A} \bar{y}_{L}+x_0, \quad \text{where} \quad 
\bar{y}_L =\argmin_{\bar{y}} \left\{ \|A L^{\dagger}_{A} \bar{y}-\bar{b}\|_2^2+{\lambda}\|\bar{y}\|_2^2\right\},
\end{equation}
where $L^{\dagger}_{A}$ is the $A$-weighted pseudoinverse of $L$, defined as $L^{\dagger}_{A} = [I - (A(I-L^{\dagger}L))^{\dagger}A]L^\dagger $ ($L^\dagger$ denotes the Moore-Penrose pseudoinverse of $L$); $x_0$ is the component of the solution $\hat{x}$ in the null space of $L^{(k)}$ and $\bar{b} = b-A x_0$. The matrix $ L^{\dagger}_{A}$ can now be considered as an (iteration dependent) right preconditioner and incorporated into the space of the solutions using flexible Krylov methods, see, e.g. \cite{CALVETTI2007378}\cite{tv_fgmres}\cite{tv_flsqr}. This strategy circumvents the need for an inner-outer scheme, and provides a much faster convergence than TV - CGSL. Moreover, in a hybrid fashion, it allows for the regularization parameter $\lambda = \lambda_k$ to be computed on-the-fly throughout the iterations. However, flexible Krylov methods require storing all the computed basis vectors so that the memory requirements increase with the number of iterations. The algorithm provided in this paper is an adaptation from \cite{tv_flsqr}, using different boundary conditions for the discrete derivative operator approximation and extending it to 3D.

\section{Numerical experiments}
In this section, three numerical experiments are presented to showcase different aspects of the Krylov subspace methods described in this paper, as well as different applications in CT. First, an example using synthetic CT data is presented to analyse the convergence and behaviour of the different algorithms: providing a comprehensive comparison including relative error and residual norm histories. Second, to validate the performance of these methods on real data, a scan of the Alderson head phantom obtained in a Philips Allura medical CT scanner is reconstructed with full sampling and under-sampled projections. Finally, a bumblebee image obtained with an industrial Nikon CT scanner is reconstructed for some of the algorithms, in a real large-scale problem. Note that the large-scale nature of the CBCT image reconstruction problems means that the implementations are often in single precision floating point arithmetic.  

The first two experiments were carried out in a laptop with a Intel Core i7-7700HQ with 16GB of RAM and a GTX 1070 NVIDIA GPU. The \textmu-CT reconstruction was performed in a machine with an AMD EPYC 7352 with 126GB of RAM and 4 NVIDIA Quadro RTX 6000.

\subsection{Comprehensive convergence comparison on synthetic data}
In this experiment we explore the behaviour of the algorithms presented in this paper in the context of 3D CBCT, using the available implementation in the TIGRE toolbox. Since this is a simulated toy example, relative residual norms, i.e. $\|Ax_i-b\|/\|b\|$, and relative error norms $\|x_i-x_{gt}\|/\|x_{gt}\|$ (for a given iteration number $i$ and a ground truth image $x_{gt}$) will be given to illustrate the behaviour of the solvers. %This can be considered as an illustrative example of the general performance of these algorithms for most CT problems.

The data presented in this example concerns the measurements of a synthetic dataset of a human head of size $64$$\times$$64$$\times$$64$. The results for this experiment are presented an analysed in three different subsections. First, the results for LSQR and hybrid LSQR (using different regularization parameter choices) are displayed to illustrate the typical behaviour of Krylov methods for CT problems, and can be observed in Figures \ref{fig:lsqr_Recon}, \ref{fig:lsqrnorm}. Second, the Krylov methods presented in this paper for the least squares problems with and without Tikhonov regularization are shown in Figures \ref{fig:algos_e0} and \ref{fig:norms}. Finally, a sparse-view version of the same simulation is used to showcase the Krylov methods that enforce TV regularization. 

\subsubsection{Illustration of Krylov methods' typical behaviour} %In this section, the results for the first experiment using LSQR and variations thereof are compared and analysed with the aim of giving an overview of Krylov methods' typical behaviour. 
\begin{figure}
    \centering
    \captionsetup[subfigure]{labelformat=empty,justification=centering}

    \subfloat[]{
    \centering
    \includegraphics[width=0.16\columnwidth,valign=c]{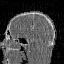}
    } \hspace*{-0.9em}
\subfloat[]{    
\centering
    \includegraphics[width=0.16\columnwidth,valign=c]{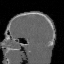}
} \hspace*{-0.9em}    
\subfloat[]{
\centering

\includegraphics[width=0.16\columnwidth,valign=c]{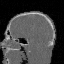}
} \hspace*{-0.9em}
\subfloat[]{
\centering
\includegraphics[width=0.16\columnwidth,valign=c]{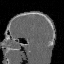}
} \hspace*{-0.9em}
\subfloat[]{
\centering
\includegraphics[width=0.16\columnwidth,valign=c]{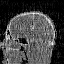}
} \hspace*{-0.9em}
\subfloat[]{
\centering
\includegraphics[width=0.16\columnwidth,valign=c]{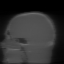}
} \hspace*{-0.9em}

\vspace{-1.7em}

    \subfloat[LSQR]{
    \centering
    \includegraphics[width=0.16\columnwidth,valign=c]{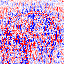}
    } \hspace*{-0.9em}
\subfloat[hybrid LSQR \footnotesize{(GCV)}]{    
\centering
    \includegraphics[width=0.16\columnwidth,valign=c]{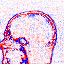}
    } \hspace*{-0.9em}
\subfloat[hybrid LSQR \footnotesize{(DP, 1.5\% noise)}]{
\centering
    \includegraphics[width=0.16\columnwidth,valign=c]{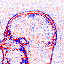}
    } \hspace*{-0.9em}
\subfloat[hybrid LSQR \footnotesize{($\lambda=20$)}]{
\centering
    \includegraphics[width=0.16\columnwidth,valign=c]{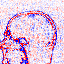}
    } \hspace*{-0.9em}
\subfloat[hybrid LSQR \footnotesize{($\lambda=2$)}]{
\centering
    \includegraphics[width=0.16\columnwidth,valign=c]{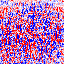}
    } \hspace*{-0.9em}
\subfloat[hybrid LSQR \footnotesize{($\lambda=200$)}]{
\centering
    \includegraphics[width=0.16\columnwidth,valign=c]{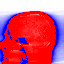}
    } \hspace*{-0.9em}
    
    \caption{Reconstruction of phantom head data using different LSQR versions and parameter selection procedures. (top row) slice of final images, shown in range [0, 1] mm\textsuperscript {-1} (bottom row) difference images w.r.t. the ground truth, shown in range [-0.1, 0.1] mm\textsuperscript {-1}.}
    \label{fig:lsqr_Recon}
\end{figure}

\begin{figure}
    \centering
    \subfloat[ ]{\includegraphics[width=0.32\columnwidth,valign=c]{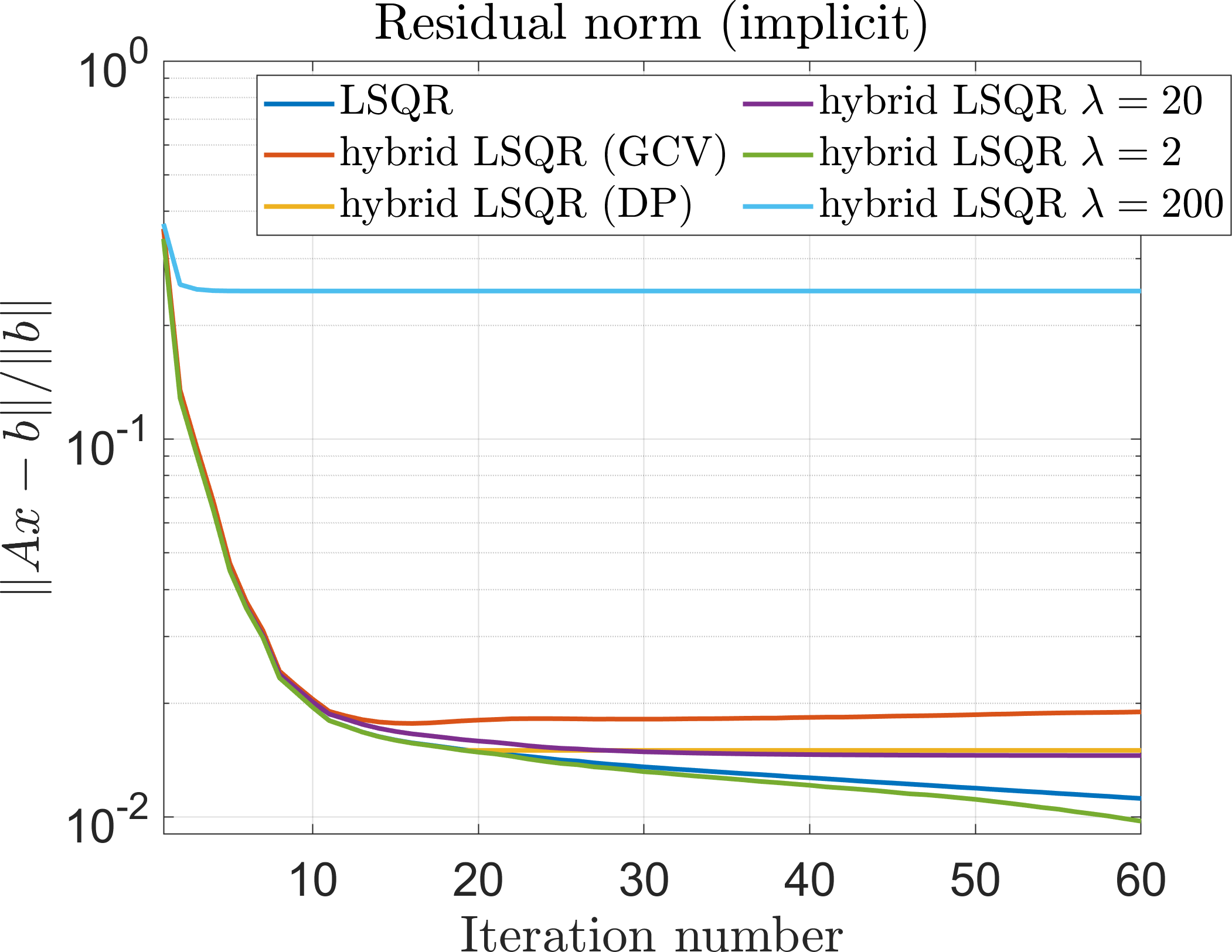}}
    ~
    \subfloat[ ]{\includegraphics[width=0.32\columnwidth,,valign=c]{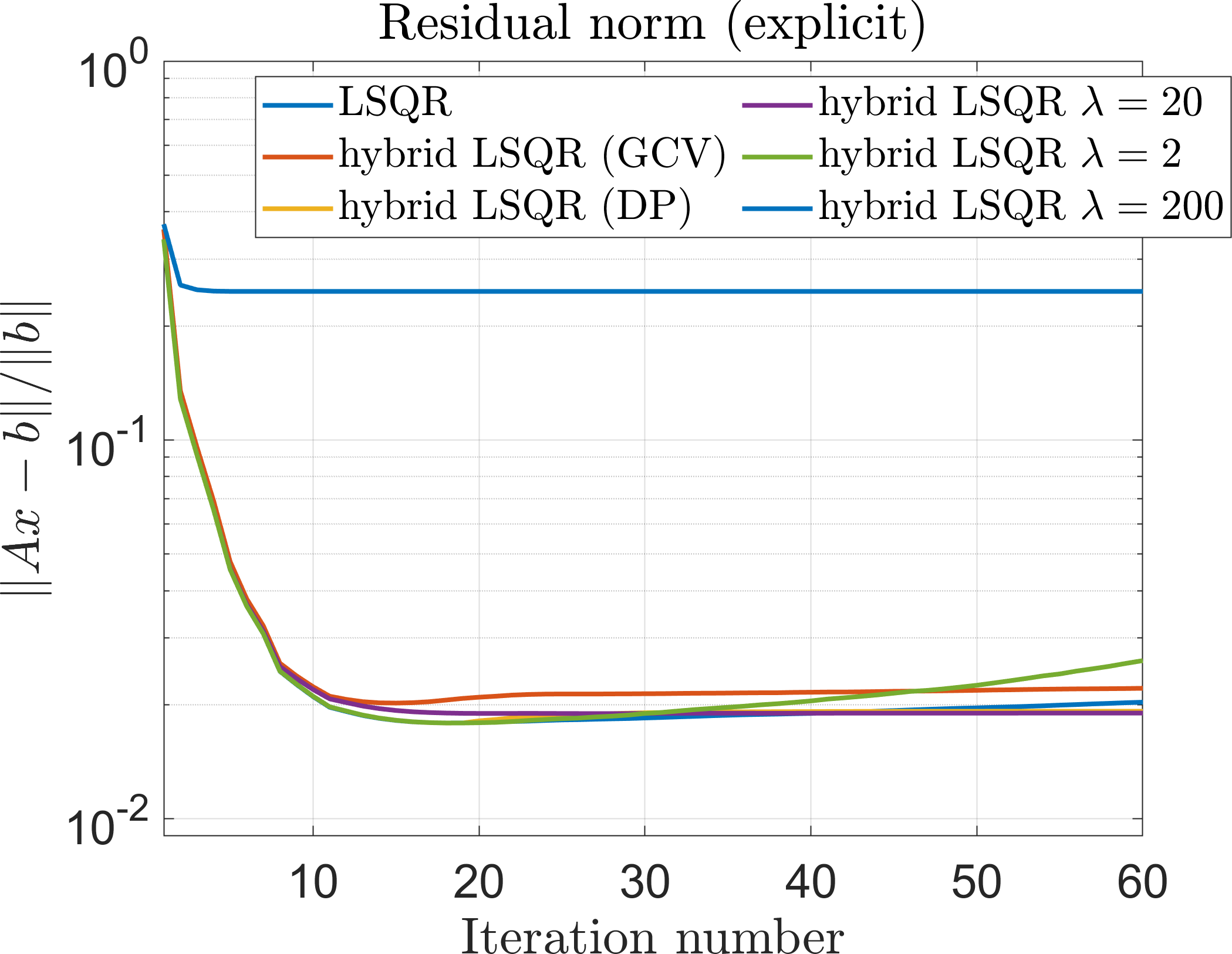}}
     ~
    \subfloat[ ]{\includegraphics[width=0.32\columnwidth,,valign=c]{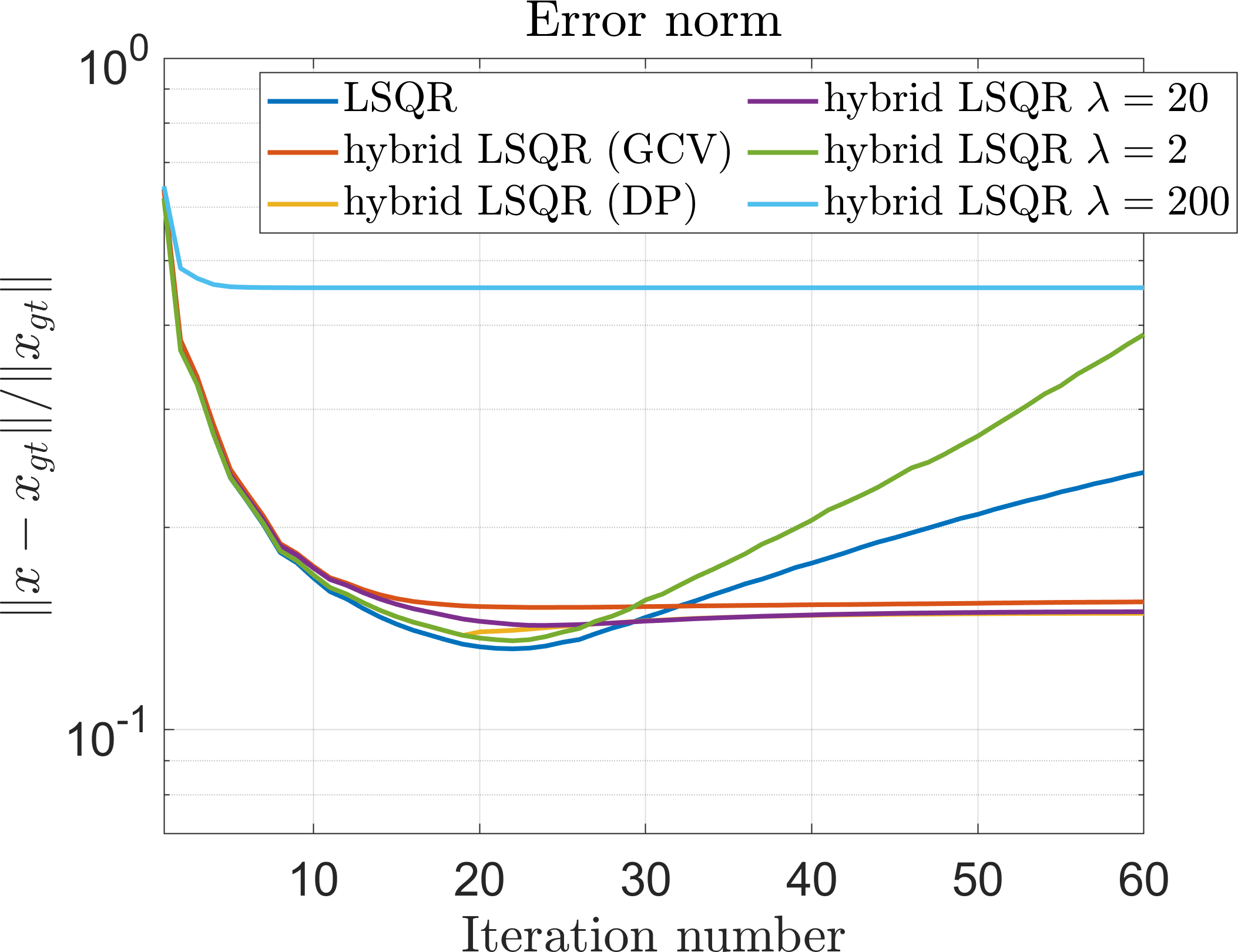}}
    \caption{(a) Implicit relative residual norms, (b) computed relative residual norms and (c) relative error norms for the algorithms of interest, per iteration.}
    \label{fig:lsqrnorm}
\end{figure}

In this experiment, 60 equidistant angles spanning the full circular range are simulated, with added Poisson noise (assuming an air photon count of $I_0=1\times10^5$) and Gaussian noise (with standard deviation of $\sigma=0.5$) which model both the photon and electronic noise expected in a CT scanner \cite{xu2009electronic}\cite{liu2012adaptive}. The reconstructions can be observed in Figure \ref{fig:lsqr_Recon}, while relative residual and error norms are displayed in Figure \ref{fig:lsqrnorm}.

First, we can observe that LSQR undergoes semiconvergence. As briefly mentioned in Section \ref{sec_lsqr}, this is characterized by the decrease of the relative error norm in the first iterations until it reaches a minima after which it starts increasing (see Figure \ref{fig:lsqrnorm}(c)); while the residual norm is in exact arithmetic minimized at each iteration and therefore decreases monotonically throughout the iterations (see Figure \ref{fig:lsqrnorm}(a)). This phenomenon is the reason why early stopping is crucial to obtain a good reconstruction in inverse problems with noisy measurements when using an iterative solver that acts directly on the least squares problem \ref{eq:lest_squares}.  

Second, we want to illustrate how the choice of a good regularization parameter is crucial to obtain a meaningful reconstruction when solving the Tikhonov problem \ref{eq:Tikhonov}. As can be observed in both the reconstructions (Figure \ref{fig:lsqr_Recon}) and the relative error norm histories (Figure \ref{fig:lsqrnorm}(c)), the semi-automatic parameter choice criteria provided in this implementation find appropriate parameters $\lambda_k$ at each iteration to obtain a good reconstruction without fine tuning. Alternatively, one can choose a parameter $\lambda$ ahead of the iterations. In this case, note that an under-regularized problem (see Figure \ref{fig:lsqr_Recon} for $\lambda=2$) will produce a noisy-looking image (in which case one should re-run the algorithms using a higher value for the parameter $\lambda$); while an over-regularised problem (see Figure \ref{fig:lsqr_Recon} for $\lambda=200$), will produce an overly smooth reconstruction (in which case one should use a smaller value for $\lambda$).

Last, a very interesting thing to note is the mismatch between the theoretical or implicit residuals in Figure \ref{fig:lsqrnorm}(a), i.e. computed using $\|\|b\|e_1-H_k y_i\|/\|b\|$ or mathematical recurrences (see \cite{lsqr} for LSQR), and the residuals computed explicitly throughout the iterations in Figure \ref{fig:lsqrnorm}(b), i.e. using $\|Ax_i-b\|/\|b\|$ directly. This can be due to loss of orthogonality (mainly attributed to the mismatched backprojector) or to an accumulation of numerical errors and precision loss (most objects are stored in single precision floating point arithmetic, with large differences in the order of magnitude of the parameters). Note that this happens both for the algorithms that incorporate re-orthogonalization and for the ones that do not. As this mismatch flags a deviation of the algorithm from its expected behaviour in exact arithmetic, TIGRE explicitly computes the residual norms at each iteration and stops the algorithm once they increase.

\subsubsection{(hybrid) Krylov methods for least squares problems}
\begin{figure}
    \centering
    \captionsetup[subfigure]{labelformat=empty,justification=centering,position=top}
    
     % Row 1
    \subfloat[SIRT]{
    \centering
    \includegraphics[width=0.19\columnwidth,valign=c]{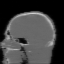}
    } \hspace*{-0.9em}
    \subfloat[CGLS]{
       \centering
    \includegraphics[width=0.19\columnwidth,valign=c]{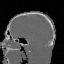}
    } \hspace*{-0.9em}
    \subfloat[LSQR]{
       \centering
    \includegraphics[width=0.19\columnwidth,valign=c]{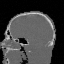}
    } \hspace*{-0.9em}
    \subfloat[LSMR \\ \footnotesize{$\lambda=0$}]{
     \centering
    \includegraphics[width=0.19\columnwidth,valign=c]{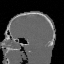}
    } \hspace*{-0.9em}
    \subfloat[LSMR \\ \footnotesize{$\lambda=30$}]{
      \centering
    \includegraphics[width=0.19\columnwidth,valign=c]{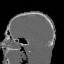}
    }%
    \vspace{-0.7em}
     \captionsetup[subfigure]{labelformat=empty,justification=centering,position=bottom}
    % Row2
    \subfloat[hybrid LSQR]{
       \centering
    \includegraphics[width=0.19\columnwidth,valign=c]{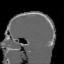}
    } \hspace*{-0.9em}
    \subfloat[AB-GMRES]{
       \centering
    \includegraphics[width=0.19\columnwidth,valign=c]{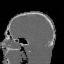}
    } \hspace*{-0.9em}
     \subfloat[BA-GMRES]{
       \centering
    \includegraphics[width=0.19\columnwidth,valign=c]{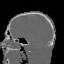}
    } \hspace*{-0.9em}
     \subfloat[AB-GRMES \\ \footnotesize{B = FDK}]{
       \centering
    \includegraphics[width=0.19\columnwidth,valign=c]{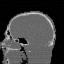}
    } \hspace*{-0.9em}
     \subfloat[BA-GRMES \\ \footnotesize{B = FDK}]{
       \centering
    \includegraphics[width=0.19\columnwidth,valign=c]{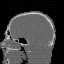}
    }
    
   \vspace{1em}
   
%% Diff:
\captionsetup[subfigure]{labelformat=empty,justification=centering,position=top}
   
     % Row 1
    \subfloat[SIRT]{
    \centering
    \includegraphics[width=0.19\columnwidth,valign=c]{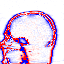}
    } \hspace*{-0.9em}
    \subfloat[CGLS]{
       \centering
    \includegraphics[width=0.19\columnwidth,valign=c]{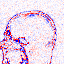}
    } \hspace*{-0.9em}
    \subfloat[LSQR]{
       \centering
    \includegraphics[width=0.19\columnwidth,valign=c]{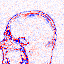}
    } \hspace*{-0.9em}
    \subfloat[LSMR \\ \footnotesize{$\lambda=0$}]{
     \centering
    \includegraphics[width=0.19\columnwidth,valign=c]{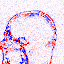}
    } \hspace*{-0.9em}
    \subfloat[LSMR \\ \footnotesize{$\lambda=30$}]{
      \centering
    \includegraphics[width=0.19\columnwidth,valign=c]{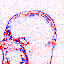}
    } \vspace{-0.7em}
     \captionsetup[subfigure]{labelformat=empty,justification=centering,position=bottom}
    % Row2
    \subfloat[hybrid LSQR]{
       \centering
    \includegraphics[width=0.19\columnwidth,valign=c]{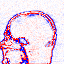}
    } \hspace*{-0.9em}
    \subfloat[AB-GMRES]{
       \centering
    \includegraphics[width=0.19\columnwidth,valign=c]{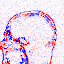}
    } \hspace*{-0.9em}
     \subfloat[BA-GMRES]{
       \centering
    \includegraphics[width=0.19\columnwidth,valign=c]{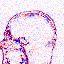}
    } \hspace*{-0.9em}
     \subfloat[AB-GRMES \\ \footnotesize{B = FDK}]{
       \centering
    \includegraphics[width=0.19\columnwidth,valign=c]{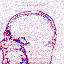}
    } \hspace*{-0.9em}
     \subfloat[BA-GRMES \\ \footnotesize{B = FDK}]{
       \centering
    \includegraphics[width=0.19\columnwidth,valign=c]{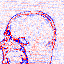}
    }
    
   \vspace{1em}

        \caption{Reconstruction of phantom head data using several Krylov methods (top) slice of final images, shown in range [0, 1] mm\textsuperscript {-1} (bottom) difference images w.r.t. the ground truth, shown in range [-0.1, 0.1] mm\textsuperscript {-1}.  }
    \label{fig:algos_e0}
\end{figure}

\begin{figure}
    \centering
    \subfloat[ ]{\includegraphics[width=0.49\columnwidth,valign=c]{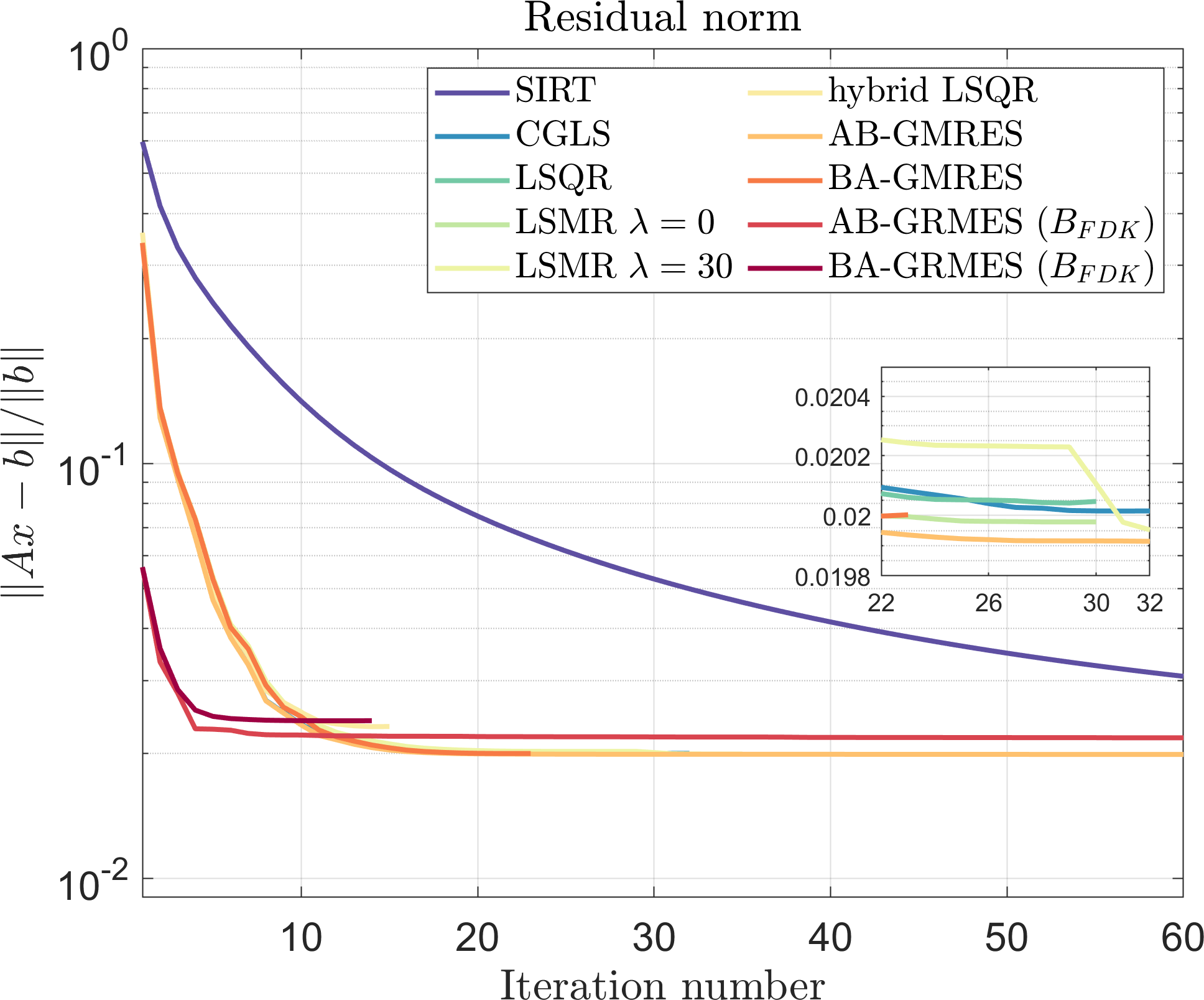}}
    ~
    \subfloat[ ]{\includegraphics[width=0.49\columnwidth,,valign=c]{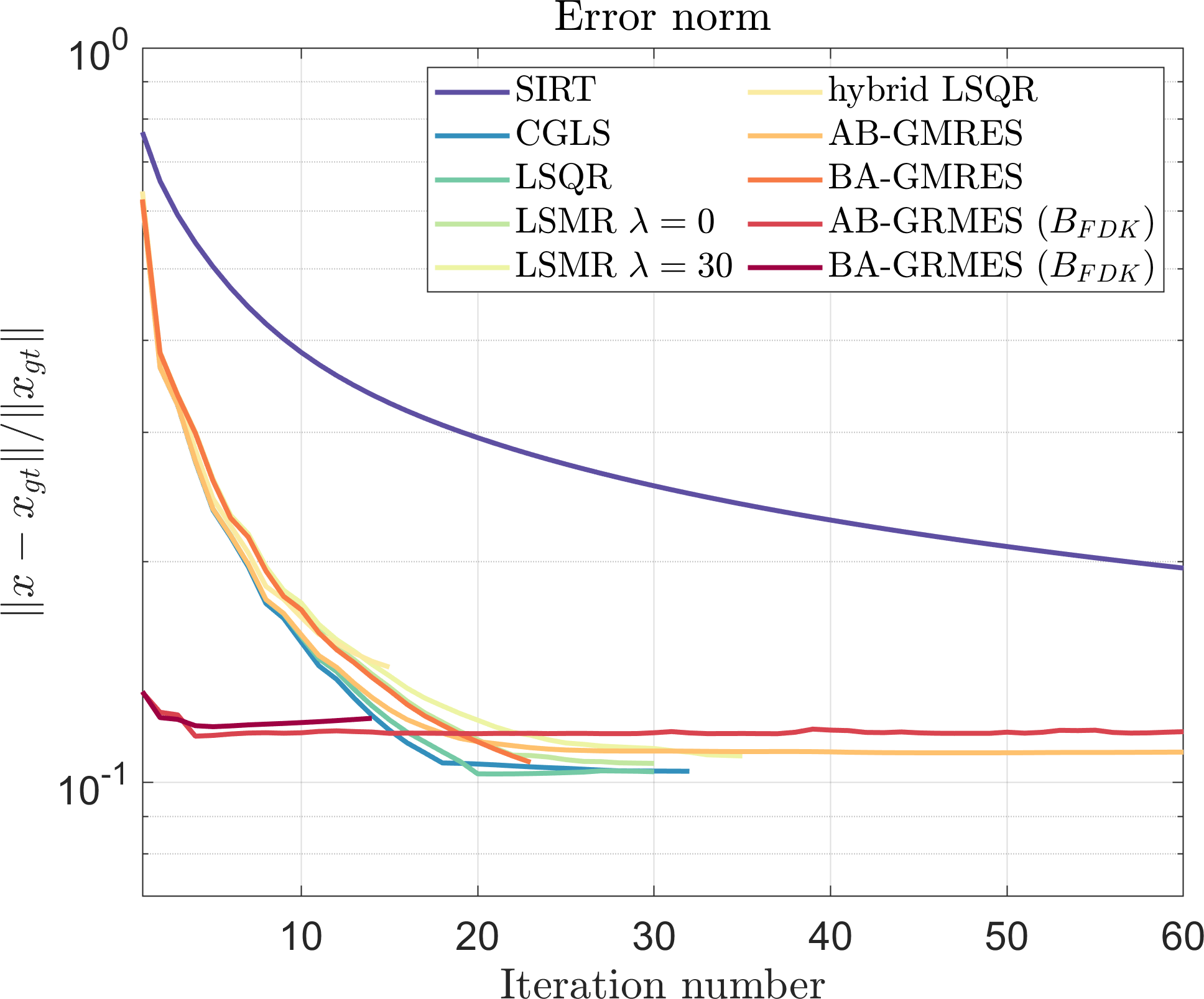}}
    \caption{(a) Relative residual norms and (b) relative error norms for the compared algorithms, per iteration.}
    \label{fig:norms}
\end{figure}

This experiment concerns a dataset with 180 equidistant angular projections with the same noise distribution used in the previous section ($I_0=1\times10^5$, $\sigma=0.5$). The results for all the algorithms presented in this work for the least squares problem with or without Tikhonov regularization, are shown for a maximum of 60 iterations. As a baseline, the results are compared to the solutions computed with SIRT: a particular choice from the most commonly used family of algorithms in CT, the SIRT-like family (SIRT, OS-SART, SART, etc) \cite{kak2001principles}, which is computationally equivalent to the Krylov methods used in this example (i.e., they have an equivalent amount of flops per iteration). Figure \ref{fig:algos_e0} shows a slice of the reconstructed image obtained using the different methods on top of its corresponding error (difference between the reconstructed slice and the ground truth). Figure \ref{fig:norms} shows the relative residual norm and relative error norm histories for all the algorithms against the iteration numbers. 

In Figure \ref{fig:norms} one can observe the very fast convergence of Krylov methods: both in terms of the relative residual norm and of the relative error norm. In this particular experiment, between 10 and 20 iterations of the compared Krylov subspace methods are sufficient to obtain a good reconstruction of the original image while, after 60 iterations, SIRT has still not converged and has failed to compute meaningful reconstruction. It can also be observed that the different Krylov methods perform similarly, with LSQR producing results of slightly better quality than CGLS in terms of error norm.

For this particular example, the iterations are stopped early if the norm of the explicit residual increases between two consecutive iterations, as this is a sign of loss of orthogonality in the basis vectors or of the accumulation of computational errors. This happens for most compared Krylov subspace algorithms, and note that for CGLS, LSMR and LSQR, as a side-effect, this leads to regularization by early stopping and avoids the semiconverge behaviour. However, for these algorithms, one should stop the iterations early even when the explicit residual norm does not increase (for example, by monitoring the stabilization of the residual norm or using other stopping criteria). It is also remarkable to observe that the AB/BA-GMRES algorithms less commonly display increasing residuals, as they alleviate the problems associated with mismatched backprojectors. 

\subsubsection{Krylov methods for total variation regularization} 
\begin{figure}
    \centering
    \captionsetup[subfigure]{labelformat=empty,justification=centering}
    
\subfloat[]{
     \centering
    \includegraphics[width=0.19\columnwidth,valign=c]{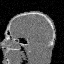}
    } \hspace*{-0.9em}
\subfloat[]{
     \centering
    \includegraphics[width=0.19\columnwidth,valign=c]{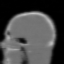}
    } \hspace*{-0.9em}
    \subfloat[]{
     \centering
    \includegraphics[width=0.19\columnwidth,valign=c]{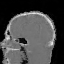}
    } 
    \vspace{-2.3em}
    \subfloat[CGLS]{
     \centering
    \includegraphics[width=0.19\columnwidth,valign=c]{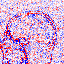}
    } \hspace*{-0.9em}
     \subfloat[CGLS-TV]{
     \centering
    \includegraphics[width=0.19\columnwidth,valign=c]{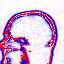}
    } \hspace*{-0.9em}
    \subfloat[hybrid fLSQR-TV]{
     \centering
    \includegraphics[width=0.19\columnwidth,valign=c]{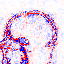}
    } 
    
    \caption{Reconstruction of phantom head data using CGLS and TV regularized Krylov methods. (top row) slice of final images, shown in range [0, 1] mm\textsuperscript {-1} (bottom row) difference images w.r.t. the ground truth, shown in range [-0.1, 0.1] mm\textsuperscript {-1}.   }
    \label{fig:algos_e3}
\end{figure}

For this experiment, the simulated CT measurements of the dataset described in the previous section are reduced and correspond to 60 equidistant projections: generating a more ill-posed problem. Moreover, the Poisson noise for this problem is increased so that $I_0=1 \times 10^4$. In this case, more prior information on the solution is needed to obtain a good reconstruction of the original image, and therefore the described methods involving TV regularization become more meaningful. In particular, CGLS, CLGS with TV regularization and hybrid fLSQR with TV regularization are showcased in this experiment. The reconstructions obtained by these methods after 60 iterations can be seen in Figure \ref{fig:algos_e3}, where the smoothing but edge preserving behaviour of the TV regularization is visibly clear. 

\begin{figure}
    \centering
    \subfloat[ ]{\includegraphics[width=0.49\columnwidth,valign=c]{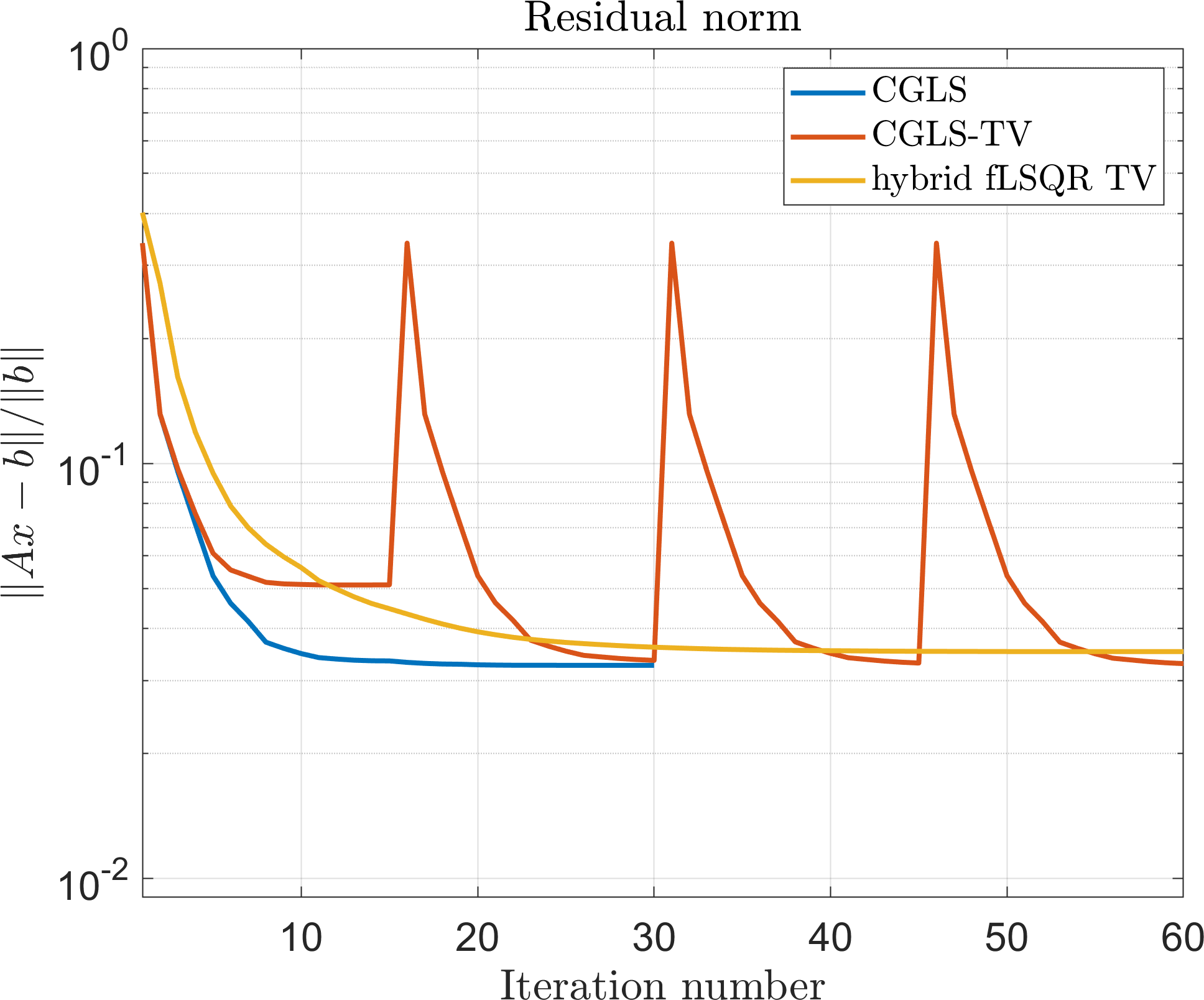}}
    ~
    \subfloat[ ]{\includegraphics[width=0.49\columnwidth,valign=c]{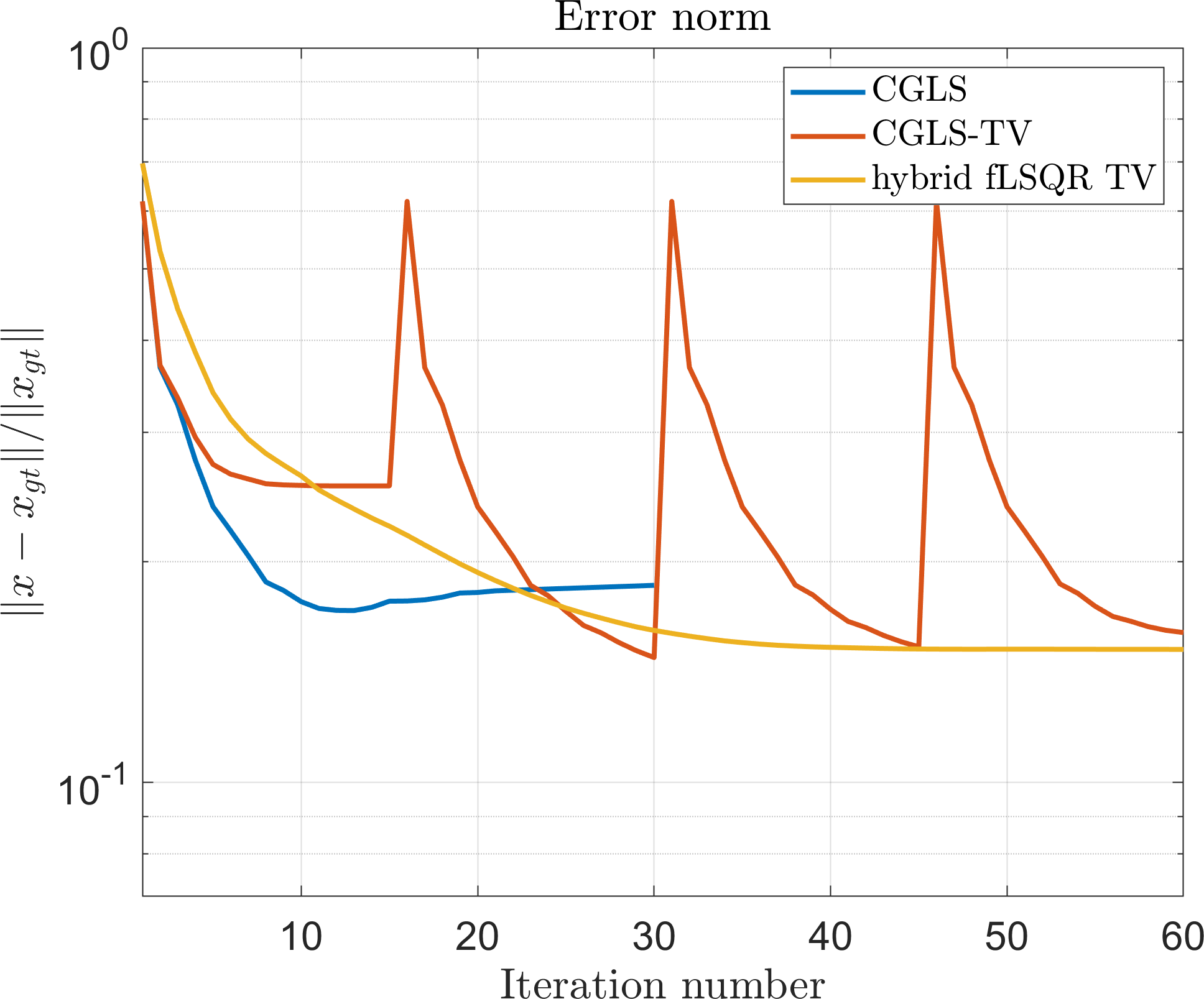}}
    \caption{(a) Relative residual norms and (b) relative error norms for the compared algorithms, per iteration.}
    \label{fig:normstv}
\end{figure}

Figure \ref{fig:normstv} shows the relative residual norms and the relative error norms throughout the iterations for the compared methods. In this example it can be clearly observed that CGLS semiconverges due to the ill-posedness of the problem and the noise in the data. It is also important to explain that the behaviour of CGLS-TV in terms or relative residual and error norms is expected. Here, the `peaks' correspond to the starts of each new cycle of inner iterations, also known as cold restarts. For this particular experiment, the number of inner iterations in each cycle is chosen a-priori to be 15 iterations, but this could also be set adaptively using a stopping criterion for the inner iterations. As long as the number of inner iterations is sufficiently large, this algorithm produces very good reconstructions with the properties expected of TV regularized solutions (this is especially desirable for highly noisy datasets). The experiment shows that likely two outer iterations (30 iterations in total) would be sufficient for this particular example. Finally, even if hybrid fLSQR-TV displays a slower decay of the residual norm, it still produces fastly decreasing error norms compared to CGLS-TV, and does not exhibit semiconvergence. Note that this method does not require to set a number of inner iterations, and the regularization parameter can be chosen semi-automatically on-the-fly. However, it requires storing all the generated basis vectors and has the additional cost of an (approximated) matrix-vector product with $L^{\dagger}_{A}(L^{\dagger}_{A})^T$ at each iteration (in the codes provided, this is done efficiently using an iterative method).

\subsection{Medical CT experiments}

This experiment concerns a medical imaging application and has the aim to highlight the performance of Krylov subspace methods on real data. In particular, the computational times for the different algorithms are given in this example to highlight the fast convergence of Krylov methods.

\begin{figure}
\centering
    % row 1
    \captionsetup[subfigure]{labelformat=empty,justification=centering,position=bottom}
    \subfloat[FDK]{%
    \centering
    \includegraphics[width=0.15\columnwidth,valign=c]{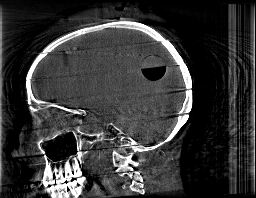}
    } \hspace*{-0.6em}
    \subfloat[SIRT \\ \tiny{30 iter}]{%
    \centering
    \includegraphics[width=0.15\columnwidth,valign=c]{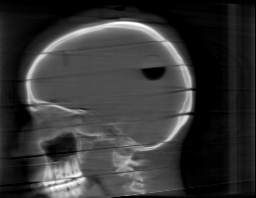}
    } \hspace*{-0.6em}
    \subfloat[SIRT\\ \tiny{150 iter}]{%
    \centering
    \includegraphics[width=0.15\columnwidth,valign=c]{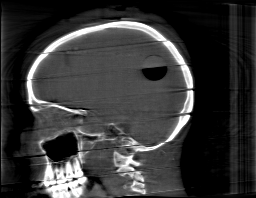}
    }\quad
    \subfloat[FDK]{%
    \centering
    \includegraphics[width=0.15\columnwidth,valign=c]{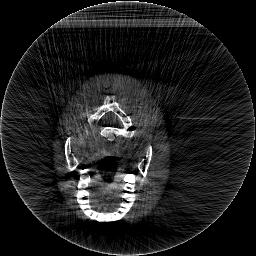}
    } \hspace*{-0.6em}
    \subfloat[SIRT \\ \tiny{30 iter}]{%
    \centering
    \includegraphics[width=0.15\columnwidth,valign=c]{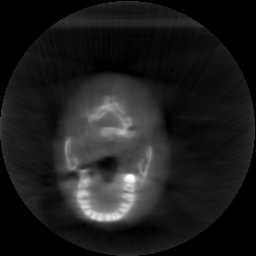}
    }\hspace*{-0.6em}
    \subfloat[SIRT\\ \tiny{150 iter}]{%
    \centering
    \includegraphics[width=0.15\columnwidth,valign=c]{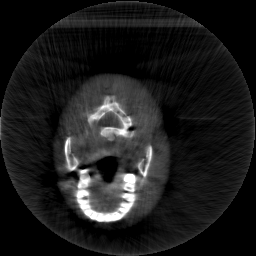}
    } %
    % row change
    \vspace{-0.9em}
    % row 2
    \captionsetup[subfigure]{labelformat=empty,justification=centering,position=bottom}
    \subfloat[OS-SART \\ \tiny{60 iter}]{%
    \centering
    \includegraphics[width=0.15\columnwidth,valign=c]{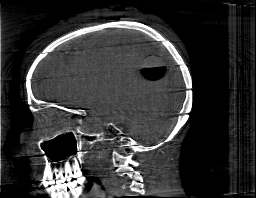}
    } \hspace*{-0.6em}
    \subfloat[CGLS \\ \tiny{30 iter}]{%
    \centering
    \includegraphics[width=0.15\columnwidth,valign=c]{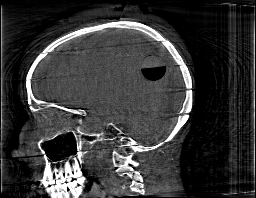}
    } \hspace*{-0.6em}
    \subfloat[LSQR\\ \tiny{30 ite}]{%
    \centering
    \includegraphics[width=0.15\columnwidth,valign=c]{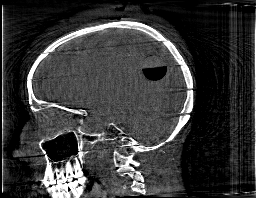}
    }\quad
    \captionsetup[subfigure]{labelformat=empty,justification=centering,position=bottom}
    \subfloat[OS-SART \\ \tiny{60 iter}]{%
    \centering
    \includegraphics[width=0.15\columnwidth,valign=c]{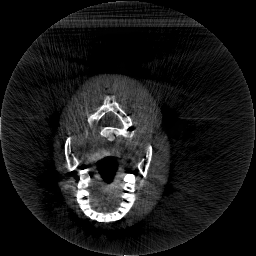}
    } \hspace*{-0.6em}
    \subfloat[CGLS \\ \tiny{30 iter}]{%
    \centering
    \includegraphics[width=0.15\columnwidth,valign=c]{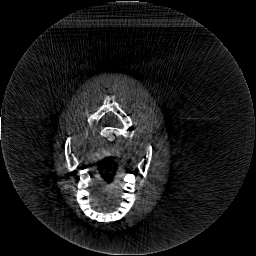}
    } \hspace*{-0.6em}
    \subfloat[LSQR\\ \tiny{30 iter}]{%
    \centering
    \includegraphics[width=0.15\columnwidth,valign=c]{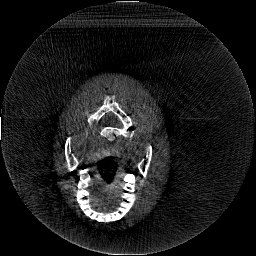}
    }%
    \vspace{-0.9em}
    % row 6
    \captionsetup[subfigure]{labelformat=empty,justification=centering,position=bottom}
    \subfloat[LSMR\\ \tiny{30 iter, $\lambda=0$}]{%
    \centering
    \includegraphics[width=0.15\columnwidth,valign=c]{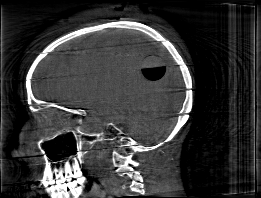}
    } \hspace*{-0.6em}
    \subfloat[LSMR\\ \tiny{30 iter, $\lambda=30$}]{%
    \centering
    \includegraphics[width=0.15\columnwidth,valign=c]{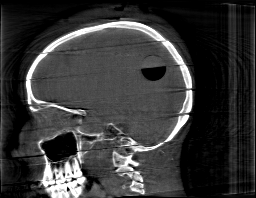}
    } \hspace*{-0.6em}
    \subfloat[hybrid LSQR \\ \tiny{30 iter}]{%
    \centering
    \includegraphics[width=0.15\columnwidth,valign=c]{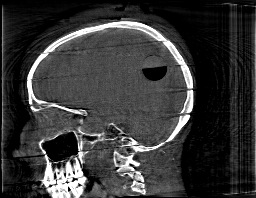}
    }\quad
       \captionsetup[subfigure]{labelformat=empty,justification=centering,position=bottom}
    \subfloat[LSMR\\ \tiny{30 iter, $\lambda=0$}]{%
    \centering
    \includegraphics[width=0.15\columnwidth,valign=c]{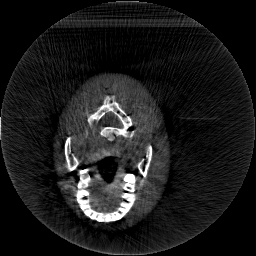}
    } \hspace*{-0.6em}
    \subfloat[LSMR\\ \tiny{30 iter, $\lambda=30$}]{%
    \centering
    \includegraphics[width=0.15\columnwidth,valign=c]{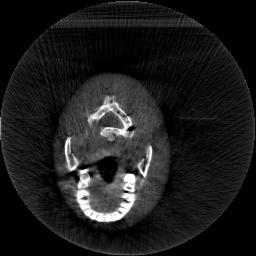}
    } \hspace*{-0.6em}
    \subfloat[hybrid LSQR \\ \tiny{30 iter}]{%
    \centering
    \includegraphics[width=0.15\columnwidth,valign=c]{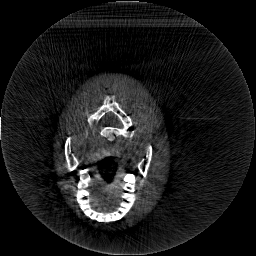}
    }
    \caption{Reconstruction of the Alderson head phantom acquired on a Phillips Allura FD20 Xper C-arm, using 289 projections. Image shown in range [0, 3] mm\textsuperscript {-1} for the (left) sagittal plane, (right) transversal plane. 
    SIRT 30 iterations and the Krylov subspace algorithms terminate within 35 seconds, while SIRT 150 iteration takes 3 minutes 50 seconds and OS-SART 6 minutes 30 seconds to converge to a solution of comparable quality.  }
    \label{fig:head_full_proj}
\end{figure}

The dataset used in this experiments consists of the Alderson head phantom measurements, acquired on a Philips Allura FD20 Xper C-arm with source settings of 80 kV and an exposure of 350 mAs, spanning a 210\textdegree{} angular range. Projections of size $512$$\times$$512$ have been used to reconstruct an image of size $256$$\times$$256$$\times$$200$ voxels. Figure \ref{fig:head_full_proj} shows two views of the image reconstructions given by different algorithms. The following explanations and comparisons are applied to both the sagittal plane (Figure \ref{fig:head_full_proj}(a)) and the transversal plane (Figure \ref{fig:head_full_proj}(b)) of the different reconstructions. The first row shows images reconstructed by FDK (considered clinical standard) and SIRT with 30 and 150 iterations, respectively. Note that the choice of 30 iterations for SIRT is taken to match the computational time required for Krylov methods to obtain a meaningful reconstruction (it can be observed that the quality of the reconstruction using SIRT in this case is not very good), while the choice of 150 iterations is taken so that SIRT displays an equivalent quality of the reconstructions than Krylov methods (taking 3 minutes and 50 seconds, almost 8 times slower than Krylov methods). In the first column of the second row, OS-SART (the ordered subset version of SIRT), is shown after 60 iterations (chosen to reach a reasonable convergence). Albeit the number of required iterations is smaller than for SIRT, OS-SART takes 6 min 30 seconds to reconstruct this image\footnote{This is specific for the particular TIGRE implementation. Faster subset algorithms can be developed using specific implementations that minimize CPU$\Leftrightarrow$GPU memory transfers. However, they require a larger amount of computations per iteration than other iterative methods, so they will still be slower than the other algorithms shown in this paper.}. In the second and third columns of the second row, the reconstructions obtained using CGLS and LSQR are shown after 30 iterations (corresponding to 30 seconds of run-time). The third row, from left to right, shows the reconstructions obtained using LSMR (with $\lambda=0$), LSMR (with $\lambda=30$) and hybrid LSQR, all of them after 30 iterations and 30 seconds of run-time. The reconstructions obtained with the studied iterative methods look less grainy than the baseline reconstruction obtained using FDK.
% {What do you mean by 'Albeit improved visibility of features is not particularly achieved, the noise level obtained by the iterative methods is visibly lower compared to FDK'.?} \ab{weird wording yes. I was just wanted to mention that the FDK image still looks kind of good, which often happens, however this is because inidvidial features (boundaries, etc) are there, and we humans are suepr good at filtering the noise. But iterative methods, while you don't see "more things", you see the things "better",i.e. with less noise. Not sur ehow to convey that without getting into an argument of visual perception psycology.}

In this experiment one can observe that iterative methods produce image reconstructions of similar or higher quality than the clinical standard FDK. Moreover, Krylov subspace methods are able to do so in significantly less computational time compared to other classic iterative reconstruction methods. This is of particular use in clinical CT, where lower reconstruction time is needed to  maximize throughput (i.e., the number of images and actions over them that can be processed per unit of time).
%, e.g., in radiation therapy reconstruction and planning often need to be done in less than 90s, or in image guided surgery, the image reconstruction should as close to  instantaneously, almost at least.

\begin{figure}
   \centering
    % row 1
    \captionsetup[subfigure]{labelformat=empty,justification=centering,position=bottom}
    \subfloat[FDK]{%
    \centering
    \includegraphics[width=0.15\columnwidth,valign=c]{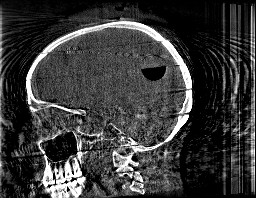}
    } \hspace*{-0.6em}
    \subfloat[SIRT \\ \tiny{30 iter}]{%
    \centering
    \includegraphics[width=0.15\columnwidth,valign=c]{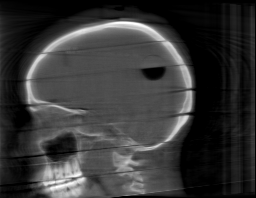}
    } \hspace*{-0.6em}
    \subfloat[SIRT\\ \tiny{150 iter}]{%
    \centering
    \includegraphics[width=0.15\columnwidth,valign=c]{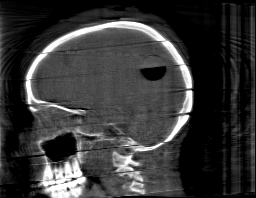}
    }\quad
    \subfloat[FDK]{%
    \centering
    \includegraphics[width=0.15\columnwidth,valign=c]{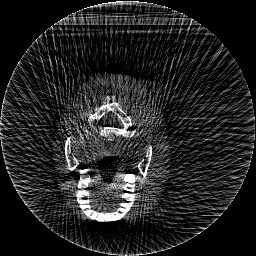}
    } \hspace*{-0.6em}
    \subfloat[SIRT \\ \tiny{30 iter}]{%
    \centering
    \includegraphics[width=0.15\columnwidth,valign=c]{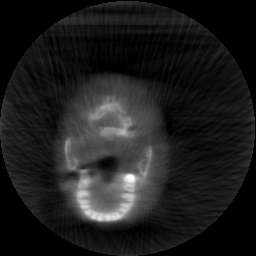}
    }\hspace*{-0.6em}
    \subfloat[SIRT\\ \tiny{150 iter}]{%
    \centering
    \includegraphics[width=0.15\columnwidth,valign=c]{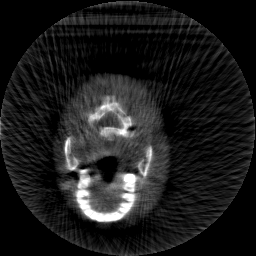}
    } %
    % row change
    \vspace{-0.9em}
    % row 2
    \captionsetup[subfigure]{labelformat=empty,justification=centering,position=bottom}
    \subfloat[OS-SART \\ \tiny{60 iter}]{%
    \centering
    \includegraphics[width=0.15\columnwidth,valign=c]{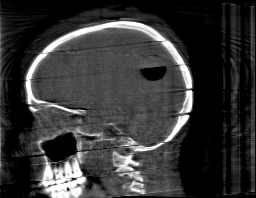}
    } \hspace*{-0.6em}
    \subfloat[CGLS \\ \tiny{30 iter}]{%
    \centering
    \includegraphics[width=0.15\columnwidth,valign=c]{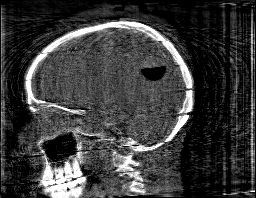}
    } \hspace*{-0.6em}
    \subfloat[LSQR\\ \tiny{30 ite}]{%
    \centering
    \includegraphics[width=0.15\columnwidth,valign=c]{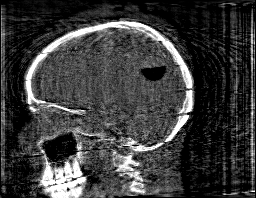}
    }\quad
    \captionsetup[subfigure]{labelformat=empty,justification=centering,position=bottom}
    \subfloat[OS-SART \\ \tiny{60 iter}]{%
    \centering
    \includegraphics[width=0.15\columnwidth,valign=c]{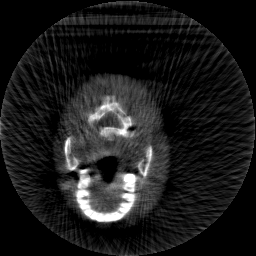}
    } \hspace*{-0.6em}
    \subfloat[CGLS \\ \tiny{30 iter}]{%
    \centering
    \includegraphics[width=0.15\columnwidth,valign=c]{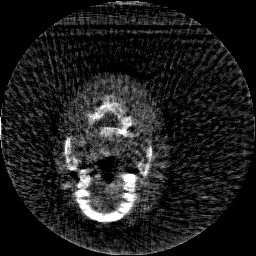}
    } \hspace*{-0.6em}
    \subfloat[LSQR\\ \tiny{30 iter}]{%
    \centering
    \includegraphics[width=0.15\columnwidth,valign=c]{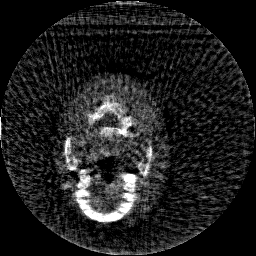}
    }%
    \vspace{-0.9em}
    % row 6
    \captionsetup[subfigure]{labelformat=empty,justification=centering,position=bottom}
    \subfloat[LSMR\\ \tiny{30 iter, $\lambda=0$}]{%
    \centering
    \includegraphics[width=0.15\columnwidth,valign=c]{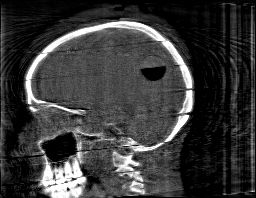}
    } \hspace*{-0.6em}
    \subfloat[LSMR\\ \tiny{30 iter, $\lambda=30$}]{%
    \centering
    \includegraphics[width=0.15\columnwidth,valign=c]{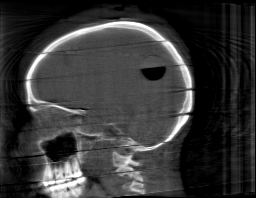}
    } \hspace*{-0.6em}
    \subfloat[hybrid LSQR \\ \tiny{30 iter}]{%
    \centering
    \includegraphics[width=0.15\columnwidth,valign=c]{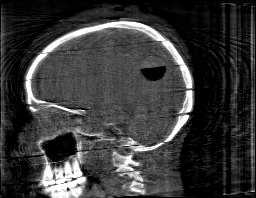}
    }\quad
       \captionsetup[subfigure]{labelformat=empty,justification=centering,position=bottom}
    \subfloat[LSMR\\ \tiny{30 iter, $\lambda=0$}]{%
    \centering
    \includegraphics[width=0.15\columnwidth,valign=c]{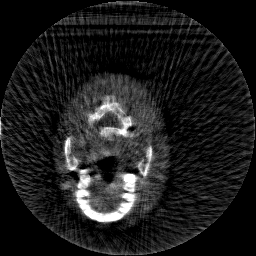}
    } \hspace*{-0.6em}
    \subfloat[LSMR\\ \tiny{30 iter, $\lambda=30$}]{%
    \centering
    \includegraphics[width=0.15\columnwidth,valign=c]{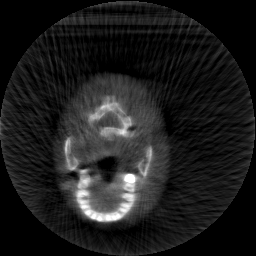}
    } \hspace*{-0.6em}
    \subfloat[hybrid LSQR \\ \tiny{30 iter}]{%
    \centering
    \includegraphics[width=0.15\columnwidth,valign=c]{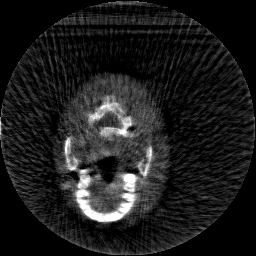}
    }
    
    \caption{Reconstruction of the Alderson head phantom acquired on a Phillips Allura FD20 Xper C-arm, using 58 projections. Image shown in range [0, 3] mm\textsuperscript {-1} for the (left) sagittal plane, (right) transversal plane.
    SIRT 30 iterations and the Krylov subspace algorithms terminate within 15 seconds, while SIRT 150 iteration takes 1 minutes 30 seconds and OS-SART 1 minutes 50 seconds to converge to a solution of comparable quality. }
    \label{fig:head_58_proj}
\end{figure}

In the following, the reconstruction results for the same example with a fifth of the projection data are shown to simulate a sparse sampling CT scan. Figure \ref{fig:head_58_proj} displays the results in the same order and for the same number of iterations already described for Figure \ref{fig:head_full_proj}. In this scenario, the Krylov subspace methods produce a good reconstruction in less than 15 seconds. Note that using SIRT in a comparable time (30 iterations) produces overly smooth reconstructions, i.e. they appear less noisy but the lack of sharpness in the edges might lead to the loss of important features in the image.  

\begin{figure}
    \centering
    \captionsetup[subfigure]{labelformat=empty,justification=centering}
    
 \subfloat[]{%
    \centering
    \includegraphics[width=0.15\columnwidth,valign=c]{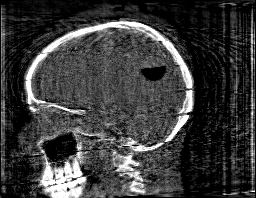}
    } \hspace*{-0.6em}
 \subfloat[]{%
    \centering
    \includegraphics[width=0.15\columnwidth,valign=c]{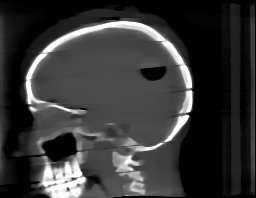}
    } \hspace*{-0.6em} \subfloat[]{%
    \centering
    \includegraphics[width=0.15\columnwidth,valign=c]{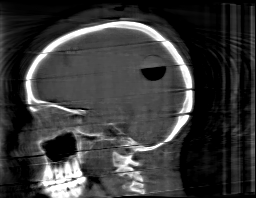}
    } 
    
    \vspace*{-2.2em}

     \subfloat[LSQR \\ \tiny{30 iter}]{%
    \centering
    \includegraphics[width=0.15\columnwidth,valign=c]{Experimet1-low_TV1.png}
    } \hspace*{-0.6em}
 \subfloat[OS-ASD-POCS \\ \tiny{60 iter}]{%
    \centering
    \includegraphics[width=0.15\columnwidth,valign=c]{Experimet1-low_TV2.png}
    } \hspace*{-0.6em} \subfloat[CGLS-TV \\ \tiny{30 iter}]{%
    \centering
    \includegraphics[width=0.15\columnwidth,valign=c]{Experimet1-low_TV3.png}
    }
    \caption{Reconstruction of the Alderson head phantom acquired on a Phillips Allura FD20 Xper C-arm, using 58 projections. Image shown in range [0, 3] mm\textsuperscript {-1} of  (top row) the sagittal plane (bottom row) the transversal plane. }
    \label{fig:head_58_proj_tv}
\end{figure}

Finally, TV regularized Krylov algorithms are used in this experiment with under-sampled projections to showcase the impact of this type of regularization on challenging CT scans. Figure \ref{fig:head_58_proj_tv} shows, for comparison, the reconstructions obtained using LSQR after 30 iterations (same as in Figure \ref{fig:head_58_proj}), OS-ASD-POCS, an ordered subset version of a well known TV regularized algorithm in tomography \cite{sidky2008image} after 60 iterations; and CGLS with TV regularization (2 outer and 15 inner iterations). For this experiment,  the running time for CGLS-TV is 1 minute, while the running time for OS-ASD-POCS is 2 minutes. It is not straightforward to establish a fair comparison purely between these two algorithms in terms of reconstruction quality, as they require the choice of different regularization hyperparameters (and there is no direct translation between the parameters for both of these algorithms). These will have a great influence in the reconstruction, balancing a closer reproduction of the fine detail features and a general smoother piece-wise constant appearance of the images. However, the results show that CGLS-TV produces good results (relatively smooth with sharp edges), while preserving the finer details of the image structures. Note that the hybrid fLSQR algorithm was not used in this experiment because the high memory needs of this algorithm were too big for the machine in which the experiments were run: this algorithm, in its current state, may not be suitable for a medical-size dataset.

\subsection{\textmu-CT scan}
\begin{figure}
    \centering
    \includegraphics[width=0.95\columnwidth]{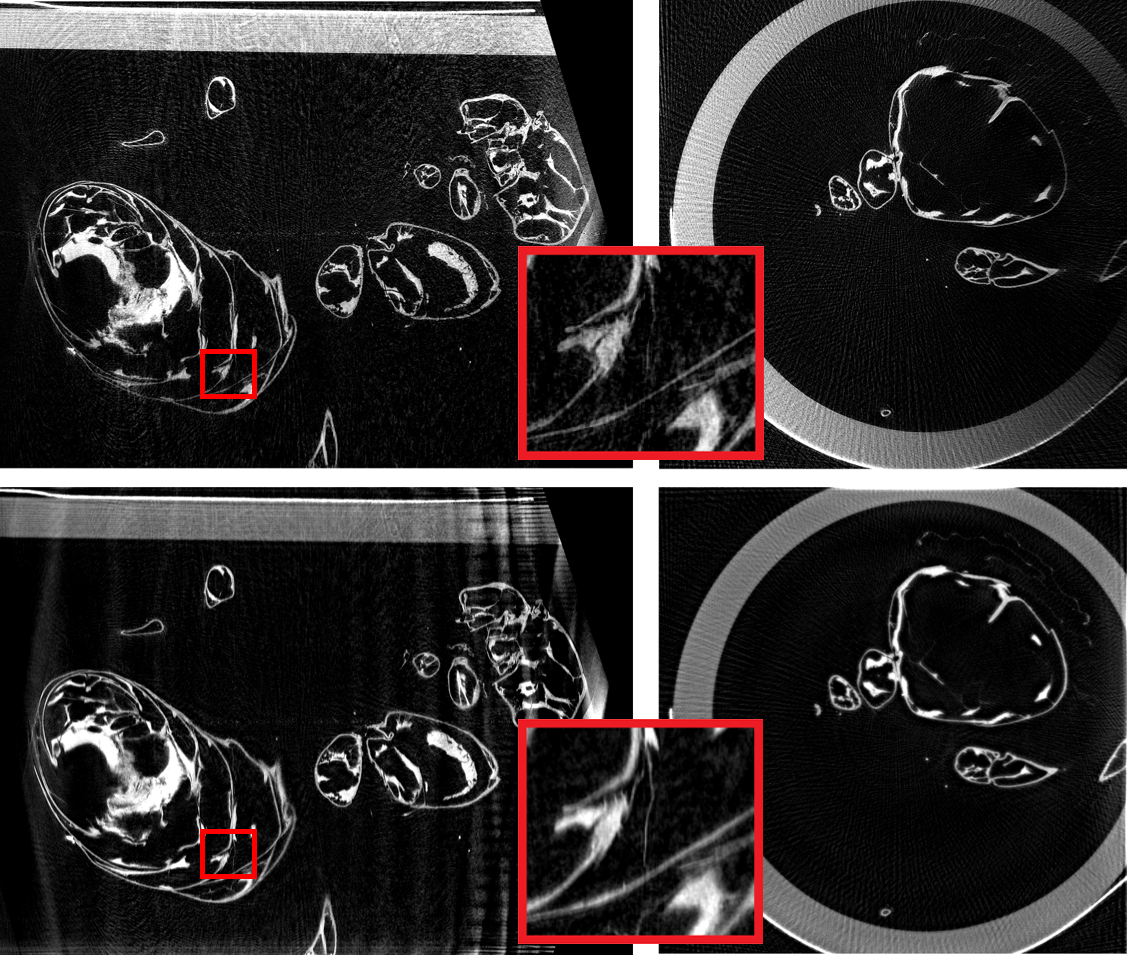}
    \caption{FDK (top) and LSQR 20 iterations (bottom) reconstruction of a \textmu-CT scan of a Bumblebee. The reconstruction is of size $1400\times1400\times2000$, using $256$ projections of a $2000\times2000$ detector. Zoomed area shows the different nature of the reconstructions, highlighting the noise rejection nature of iterative reconstruction. The small line in the center of the zoomed area can be better distinguished in the LSQR image. Image displayed at [0, 0.1] mm\textsuperscript{-1}.}
    \label{fig:bee}
\end{figure}

This experiment showcases the use of the methods presented in this work in very large scale problems where the radiation dose is not an issue. The scanned object is a wild buff-tailed bumblebee (\textit{bombus terrestris})\footnote{The unfortunate individual was found dead in the X-ray CT laboratory after getting trapped inside and became an essential part of the laboratory as an independent research dataset.}, scanned on a Nikon HMX 225 kVp CT scanner at 40 kVp with a molybdenum target. The detector was a Perkin Elmer 1621 with a gadolinium oxysulphide scintillator. The detector is of size $2000\times2000$ and we used 256 projections uniformly distributed around the circle. The reconstructed image is $1400\times1400\times2000$ with a resolution of 11.8 \textmu m per voxel. Figure \ref{fig:bee} shows the FDK reconstruction and LSQR reconstruction with 20 iterations (15 minutes of computational time). The different nature of the reconstructed image can be seen. In particular, the attenuation values of the tissue of the bumblebee are more uniform in LSQR (uniformity in the tissues is the expected result) and some features are better distinguished from the noise, particularly noticeable in the middle thin string-like structure in the zoomed-in area. 
 
\section{Discussion}
 {This section provides a discussion on some aspects reported in the numerical experiments, as well as some guidance on how to use some of these methods, with the aim of explaining potential issues that one might encounter when applying these algorithms to other datasets.}

One of the results that is mentioned in this paper is the fact that, in practice, the real residual norm can increase throughout the iterations due to a loss of orthogonality (this is not expected in exact arithmetic for the methods that theoretically minimize the residual at each iteration) or due to an accumulation numerical errors (this behaviour is also observed for the algorithms that incorporate re-orthogonalization). As described in the previous sections, this is mostly due to the use of an unmatched backprojector. The TIGRE toolbox provides an approximation of a matched backprojector \cite{matched} that mitigates the diverging behaviour in the implicit and explicit residual norms for the Krylov methods. Similarly, this is also mitigated when using algorithms that incorporate re-orthogonalization (but these come with the added cost of having to store all the computed basis vectors). An even better solution would involve implementing matched projection/backprojection operators, such as the distance driven projectors \cite{de2004distance}, or pixel driven matched projector approximations \cite{huberpixel}. Further research on the impact of numerical precision on Krylov methods would also be beneficial.

A natural question that can arise from applied scientist is which algorithm is the ``best''. First, this is an unanswerable question in general, as the algorithm choice (as well as the desired type and level of regularization) should heavily depend on the specific problem in mind and the purpose of the reconstruction. The objective of this work is not to provide such answers, but to supply easy to use and reproducible tools for exploration. However, some tips can be provided on the general use of these algorithms. It is recommended to use LSQR over CGLS, as it is a more stable algorithm but they are mathematically equivalent. In general, for 
severely undersampled CT measurements, especially with high noise levels, explicit regularization is recommended. In particular, TV regularization can be very beneficial to promote sharp edges (note that this is not only true for Krylov methods). Moreover, high regularization will produce less grainy reconstructions, but over-regularizing might lead to loosing tissue/material texture properties; this can be beneficial in some contexts, e.g. segmentation or classification, but a problem in other contexts, e.g. when the finer details are important and the noisy appearance of the reconstruction is not a problem. Finally, if enough memory is available (an image per iteration), algorithms with explicit re-orthogonalization are recommended, such as AB/BA-GMRES, to mitigate the problems derived from the loss of orthogonality.
% converge to the true solution, unlike the others\footnote{Extensive numerical results suggest that a sufficiently close enough solution is achieved if the algorithm is stopped early, or when the explicit computed residual increases}. 

It is also important to state that this is a representative but by all means not comprehensive list of all Krylov methods and their corresponding features; such as stopping criteria, see, e.g. \cite{9732394} or parameter selection criteria, see e.g. \cite{https://doi.org/10.48550/arxiv.2105.07221}\cite{10.1002/gamm.202000017}. The authors encourage the public to contribute to this work by submitting new algorithms or improving the ones available at the TIGRE toolbox.

\section{Conclusions}
 {This work describes and compares a variety of Krylov subspace methods for applied large-scale 3D CT and CBCT reconstruction, some of them used in this context for the first time. In particular, the methods included in this work are summarized in Table \ref{table:solvers}.} \\

 {The considered Krylov methods are compared and discussed in the numerical experiments, see Section 3, where it can be clearly observed that the main strength of Krylov subspace methods is their fast convergence compared to the most commonly used SIRT-like methods in iterative CT reconstruction. This is of crucial importance in medical applications, for example in image guided therapies where almost real time reconstructions are needed, but also in industrial applications where a high number of iterations is unfeasible due to the big dimensionality of the problems.
\\}

 {Finally, all the results shown in the paper are reproducible, and all the methods are provided as open source and freely accessible algorithms within the framework of the TIGRE toolbox. Some guidance on how to use these methods and a small discussion on potential results and issues one might encounter when using them on other datasets is given in the discussion, see Section 4. All the methods presented in this work can be found at \url{github.com/CERN/TIGRE} under a permissive BSD-3 clause license.}

\begin{table}[h]
\renewcommand{\arraystretch}{1.5}
\small
\caption{\label{table:solvers} { List of iterative methods detailed in the paper. Here `Objective' describes the optimization problem that is solved: LS referring to the least-squares problem \ref{eq:lest_squares}, Tikh. referring to the Tikhonov-regularized least-squares problem \ref{eq:Tikhonov}, hybrid referring to adding Tikhonov regularization to the projected problem at each iteration \ref{eq:hybrid_LSQR}, TV referring to the least-squares problems with added total variation regularization \ref{eq:TV_reweighted}. Note that Tikhonov regularization for a known regularization parameter $\lambda$ can also be applied considering the augmented system \ref{eq:Tikhonov_augmented} and using any solver for LS. Finally, note that algorithms cited with an asterisk required a significant adaptation from the original papers.}}
\begin{tabular}{|p{23mm}|p{88mm}|p{21mm}|p{7mm}|}\hline
\textbf{Method} & \textbf{Description} & \textbf{Objective} & \textbf{Ref.} \\ \hline\hline
CGLS & Conjugate gradient method applied to
 the normal equations. Minimizes the residual norm. & LS & \cite{Hestenes1952MethodsOC} \\
LSQR &  Mathematically equivalent to CGLS, using GK bidiagonalization, implemented with short recursions. Minimizes the residual norm. &  LS & \cite{lsqr} \\
LSMR & Algortihm based on the GK bidiagonalization, minimizes the normal equations residual norm. The regularization parameter $\lambda$ can be provided ahead of the iterations. & LS ($\lambda$$=$$0$) \newline Tikh.\ ($\lambda$$\neq$$0$) &  \cite{lsmr} \\
AB-GMRES \newline BA-GMRES & Adaptations of GMRES (minimal residual method using Arnoldi decomposition) using a given approximation of the backprojector as either left or right preconditioning. It is more robust for unmatched backprojectors. & LS  & \cite{HANSEN2022114352} \\
hybrid LSQR &  Hybrid version of LSQR to solve Tikhonov regularized problems. The regularization parameter $\lambda$ can be chosen ahead of the iterations or using a param. choice criteria (DP or GCV).
  &   hybrid & \cite{10.1145/355993.356000} \\
TV-CGLS & Approximation of TV using a sequence of quadratic tangent majorants that are solved with CGLS.  & TV & \cite{IRN}* \\
TV-FLSQR & Approximation of TV using a sequence of quadratic tangent majorants that are partially solved throughout the iterations using FLSQR. It is faster than TV-CGLS but has a high storage cost. & TV & \cite{tv_flsqr}* \\
\hline
\end{tabular}
\end{table}
\section*{Acknowledgements}

MSL gratefully acknowledges support from the CMIH, University of Cambridge. AB acknowledges the support of EPSRC grant EP/W004445/1. CBS acknowledges support from the Philip Leverhulme Prize, the Royal Society Wolfson Fellowship, the EPSRC advanced career fellowship EP/V029428/1, EPSRC grants EP/S026045/1 and EP/T003553/1, EP/N014588/1, EP/T017961/1, the Wellcome Innovator Awards 215733/Z/19/Z and 221633/Z/20/Z, the European Union Horizon 2020 research and innovation programme under the Marie Skodowska-Curie grant agreement No. 777826 NoMADS, the Cantab Capital Institute for the Mathematics of Information and the Alan Turing Institute. The support from the personnel of the Institute of Diagnostic and Interventional Radiology and Nuclear Medicine, Wiener Neustadt, Austria, for the performance of measurements for the Alderson head phantom is gratefully appreciated. 
\FloatBarrier

\section*{Bibliography} 

\bibliographystyle{unsrt}
\bibliography{main}

\end{document}